\documentclass[11pt]{amsart}
\usepackage{amsmath}
\usepackage{amssymb}
\usepackage{amscd}
\usepackage{xspace}
\usepackage{verbatim}

\begin{document}
\numberwithin{equation}{section}
\title[Boundary trace in dihedral domains]{Boundary trace of positive solutions \\  of supercritical semilinear
elliptic\\  equations in dihedral domains}
\author{Moshe Marcus}
\address{Department of Mathematics, Technion
 Haifa, ISRAEL}
 \email{marcusm@math.technion.ac.il}
\author{Laurent Veron}
\address{Laboratoire de Math\'ematiques, Facult\'e des Sciences
Parc de Grandmont, 37200 Tours, FRANCE}
\email{veronl@lmpt.univ-tours.fr}
\thanks{Both authors were partially sponsored by the French -- Israeli cooperation program through grant
No.~3-4299.
 The first author (MM) also
wishes to acknowledge the support of the Israeli Science Foundation
through grant No.~145-05. }

\date{}

\newcommand{\txt}[1]{\;\text{ #1 }\;}
\newcommand{\tbf}{\textbf}
\newcommand{\tit}{\textit}
\newcommand{\tsc}{\textsc}
\newcommand{\trm}{\textrm}
\newcommand{\mbf}{\mathbf}
\newcommand{\mrm}{\mathrm}
\newcommand{\bsym}{\boldsymbol}
\newcommand{\scs}{\scriptstyle}
\newcommand{\sss}{\scriptscriptstyle}
\newcommand{\txts}{\textstyle}
\newcommand{\dsps}{\displaystyle}
\newcommand{\fnz}{\footnotesize}
\newcommand{\scz}{\scriptsize}
\newcommand{\be}{
\begin{equation}
}
\newcommand{\bel}[1]{
\begin{equation}
\label{#1}}
\newcommand{\ee}{
\end{equation}
}
\newcommand{\eqnl}[2]{
\begin{equation}
\label{#1}{#2}
\end{equation}
}
\newtheorem{subn}{\name}
\renewcommand{\thesubn}{}
\newcommand{\bsn}[1]{\def\name{#1}
\begin{subn}}
\newcommand{\esn}{
\end{subn}}
\newtheorem{sub}{\name}[section]
\newcommand{\dn}[1]{\def\name{#1}}   
\newcommand{\bs}{
\begin{sub}}
\newcommand{\es}{
\end{sub}}
\newcommand{\bsl}[1]{
\begin{sub}\label{#1}}
\newcommand{\bth}[1]{\def\name{Theorem}
\begin{sub}\label{t:#1}}
\newcommand{\blemma}[1]{\def\name{Lemma}
\begin{sub}\label{l:#1}}
\newcommand{\bcor}[1]{\def\name{Corollary}
\begin{sub}\label{c:#1}}
\newcommand{\bdef}[1]{\def\name{Definition}
\begin{sub}\label{d:#1}}
\newcommand{\bprop}[1]{\def\name{Proposition}
\begin{sub}\label{p:#1}}
\newcommand{\R}{\eqref}
\newcommand{\rth}[1]{Theorem~\ref{t:#1}}
\newcommand{\rlemma}[1]{Lemma~\ref{l:#1}}
\newcommand{\rcor}[1]{Corollary~\ref{c:#1}}
\newcommand{\rdef}[1]{Definition~\ref{d:#1}}
\newcommand{\rprop}[1]{Proposition~\ref{p:#1}}
\newcommand{\BA}{
\begin{array}}
\newcommand{\EA}{
\end{array}}
\newcommand{\BAN}{\renewcommand{\arraystretch}{1.2}
\setlength{\arraycolsep}{2pt}
\begin{array}}
\newcommand{\BAV}[2]{\renewcommand{\arraystretch}{#1}
\setlength{\arraycolsep}{#2}
\begin{array}}
\newcommand{\BSA}{
\begin{subarray}}
\newcommand{\ESA}{\end{subarray}}
\newcommand{\BAL}{\begin{aligned}}
\newcommand{\EAL}{\end{aligned}}
\newcommand{\BALG}{\begin{alignat}}
\newcommand{\EALG}{\end{alignat}}
\newcommand{\BALGN}{\begin{alignat*}}
\newcommand{\EALGN}{\end{alignat*}}
\newcommand{\qeda}{\hspace{10mm}\hfill $\square$}
\newcommand{\Remark}{\note{Remark}}
\newcommand{\forevery}{\quad \forall}
\newcommand{\set}[1]{\{#1\}}
\newcommand{\lra}{\longrightarrow}
\newcommand{\sgn}{\rm{sgn}}
\newcommand{\lla}{\longleftarrow}
\newcommand{\llra}{\longleftrightarrow}
\newcommand{\Lra}{\Longrightarrow}
\newcommand{\Lla}{\Longleftarrow}
\newcommand{\Llra}{\Longleftrightarrow}
\newcommand{\warrow}{\rightharpoonup}
\newcommand{
\paran}[1]{\left (#1 \right )}
\newcommand{\sqbr}[1]{\left [#1 \right ]}
\newcommand{\curlybr}[1]{\left \{#1 \right \}}
\newcommand{\abs}[1]{\left |#1\right |}
\newcommand{\norm}[1]{\left \|#1\right \|}
\newcommand{\paranb}[1]{\big (#1 \big )}
\newcommand{\lsqbrb}[1]{\big [#1 \big ]}
\newcommand{\lcurlybrb}[1]{\big \{#1 \big \}}
\newcommand{\absb}[1]{\big |#1\big |}
\newcommand{\normb}[1]{\big \|#1\big \|}
\newcommand{
\paranB}[1]{\Big (#1 \Big )}
\newcommand{\absB}[1]{\Big |#1\Big |}
\newcommand{\normB}[1]{\Big \|#1\Big \|}

\newcommand{\thkl}{\rule[-.5mm]{.3mm}{3mm}}
\newcommand{\thknorm}[1]{\thkl #1 \thkl\,}
\newcommand{\trinorm}[1]{|\!|\!| #1 |\!|\!|\,}
\newcommand{\bang}[1]{\langle #1 \rangle}
\def\angb<#1>{\langle #1 \rangle}
\newcommand{\vstrut}[1]{\rule{0mm}{#1}}
\newcommand{\rec}[1]{\frac{1}{#1}}
\newcommand{\opname}[1]{\mbox{\rm #1}\,}
\newcommand{\supp}{\opname{supp}}
\newcommand{\dist}{\opname{dist}}
\newcommand{\myfrac}[2]{{\displaystyle \frac{#1}{#2} }}
\newcommand{\myint}[2]{{\displaystyle \int_{#1}^{#2}}}
\newcommand{\mysum}[2]{{\displaystyle \sum_{#1}^{#2}}}
\newcommand {\dint}{{\displaystyle \int\!\!\int}}
\newcommand{\q}{\quad}
\newcommand{\qq}{\qquad}
\newcommand{\hsp}[1]{\hspace{#1mm}}
\newcommand{\vsp}[1]{\vspace{#1mm}}
\newcommand{\prt}{\partial}
\newcommand{\sms}{\setminus}
\newcommand{\ems}{\emptyset}
\newcommand{\ti}{\times}
\newcommand{\nind}{\noindent}
\newcommand{\pr}{^\prime}
\newcommand{\ppr}{^{\prime\prime}}
\newcommand{\tl}{\tilde}
\newcommand{\wtl}{\widetilde}
\newcommand{\sbs}{\subset}
\newcommand{\sbeq}{\subseteq}
\newcommand{\indx}[1]{_{\scriptscriptstyle #1}}
\newcommand{\ovl}{\overline}
\newcommand{\unl}{\underline}
\newcommand{\nin}{\not\in}
\newcommand{\pfrac}[2]{\genfrac{(}{)}{}{}{#1}{#2}}

\def\ga{\alpha}     \def\gb{\beta}       \def\gg{\gamma}
\def\gc{\chi}       \def\gd{\delta}      \def\ge{\epsilon}
\def\gth{\theta}                         \def\vge{\varepsilon}
\def\gf{\phi}       \def\vgf{\varphi}    \def\gh{\eta}
\def\gi{\iota}      \def\gk{\kappa}      \def\gl{\lambda}
\def\gm{\mu}        \def\gn{\nu}         \def\gp{\pi}
\def\vgp{\varpi}    \def\gr{\rho}        \def\vgr{\varrho}
\def\gs{\sigma}     \def\vgs{\varsigma}  \def\gt{\tau}
\def\gu{\upsilon}   \def\gv{\vartheta}   \def\gw{\omega}
\def\gx{\xi}        \def\gy{\psi}        \def\gz{\zeta}
\def\Gg{\Gamma}     \def\Gd{\Delta}      \def\Gf{\Phi}
\def\Gth{\Theta}
\def\Gl{\Lambda}    \def\Gs{\Sigma}      \def\Gp{\Pi}
\def\Gw{\Omega}     \def\Gx{\Xi}         \def\Gy{\Psi}

\def\CS{{\mathcal S}}   \def\CM{{\mathcal M}}   \def\CN{{\mathcal N}}
\def\CR{{\mathcal R}}   \def\CO{{\mathcal O}}   \def\CP{{\mathcal P}}
\def\CA{{\mathcal A}}   \def\CB{{\mathcal B}}   \def\CC{{\mathcal C}}
\def\CD{{\mathcal D}}   \def\CE{{\mathcal E}}   \def\CF{{\mathcal F}}
\def\CG{{\mathcal G}}   \def\CH{{\mathcal H}}   \def\CI{{\mathcal I}}
\def\CJ{{\mathcal J}}   \def\CK{{\mathcal K}}   \def\CL{{\mathcal L}}
\def\CT{{\mathcal T}}   \def\CU{{\mathcal U}}   \def\CV{{\mathcal V}}
\def\CZ{{\mathcal Z}}   \def\CX{{\mathcal X}}   \def\CY{{\mathcal Y}}
\def\CW{{\mathcal W}} \def\CQ{{\mathcal Q}}
\def\BBA {\mathbb A}   \def\BBb {\mathbb B}    \def\BBC {\mathbb C}
\def\BBD {\mathbb D}   \def\BBE {\mathbb E}    \def\BBF {\mathbb F}
\def\BBG {\mathbb G}   \def\BBH {\mathbb H}    \def\BBI {\mathbb I}
\def\BBJ {\mathbb J}   \def\BBK {\mathbb K}    \def\BBL {\mathbb L}
\def\BBM {\mathbb M}   \def\BBN {\mathbb N}    \def\BBO {\mathbb O}
\def\BBP {\mathbb P}   \def\BBR {\mathbb R}    \def\BBS {\mathbb S}
\def\BBT {\mathbb T}   \def\BBU {\mathbb U}    \def\BBV {\mathbb V}
\def\BBW {\mathbb W}   \def\BBX {\mathbb X}    \def\BBY {\mathbb Y}
\def\BBZ {\mathbb Z}

\def\GTA {\mathfrak A}   \def\GTB {\mathfrak B}    \def\GTC {\mathfrak C}
\def\GTD {\mathfrak D}   \def\GTE {\mathfrak E}    \def\GTF {\mathfrak F}
\def\GTG {\mathfrak G}   \def\GTH {\mathfrak H}    \def\GTI {\mathfrak I}
\def\GTJ {\mathfrak J}   \def\GTK {\mathfrak K}    \def\GTL {\mathfrak L}
\def\GTM {\mathfrak M}   \def\GTN {\mathfrak N}    \def\GTO {\mathfrak O}
\def\GTP {\mathfrak P}   \def\GTR {\mathfrak R}    \def\GTS {\mathfrak S}
\def\GTT {\mathfrak T}   \def\GTU {\mathfrak U}    \def\GTV {\mathfrak V}
\def\GTW {\mathfrak W}   \def\GTX {\mathfrak X}    \def\GTY {\mathfrak Y}
\def\GTZ {\mathfrak Z}   \def\GTQ {\mathfrak Q}

\font\Sym= msam10 
\def\SYM#1{\hbox{\Sym #1}}
\newcommand{\tin}{\to\infty}
\newcommand{\ssub}[1]{_{_{\! #1}}}
\newcommand{\chr}[1]{\chi\indx{#1}}
\newcommand{\rest}[1]{\big |\indx{#1}}
\newcommand{\bdw}{\prt\Gw\xspace}
\newcommand{\wkc}{weak convergence\xspace}
\newcommand{\wrto}{with respect to\xspace}
\newcommand{\cons}{consequence\xspace}
\newcommand{\consy}{consequently\xspace}
\newcommand{\Consy}{Consequently\xspace}
\newcommand{\Essy}{Essentially\xspace}
\newcommand{\essy}{essentially\xspace}
\newcommand{\mnz}{minimizer\xspace}
\newcommand{\sth}{such that\xspace}
\newcommand{\ngh}{neighborhood\xspace}
\newcommand{\nghs}{neighborhoods\xspace}
\newcommand{\seq}{sequence\xspace}
\newcommand{\seqs}{sequences\xspace}
\newcommand{\sseq}{subsequence\xspace}
\newcommand{\ifif}{if and only if\xspace}
\newcommand{\suff}{sufficiently\xspace}
\newcommand{\abc}{absolutely continuous\xspace}
\newcommand{\sol}{solution\xspace}
\newcommand{\subss}{sub-solutions\xspace}
\newcommand{\subs}{sub-solution\xspace}
\newcommand{\supers}{super-solution\xspace}
\newcommand{\superss}{super-solutions  \xspace}
\newcommand{\Wlg}{Without loss of generality\xspace}
\newcommand{\wlg}{without loss of generality\xspace}
\newcommand{\locun}{locally uniformly\xspace}
\newcommand{\bvp}{boundary value problem\xspace}
\newcommand{\bvps}{boundary value problems\xspace}
\def\RN{\BBR^N}
\def\({{\rm (}}
\def\){{\rm )}}
\def\loc{\indx{\rm loc}}
\def\bmn{\mathbf{n}}
\def\bma{\mathbf{a}}
\def\prtn{\prt_{\bmn}}
\def\1{\\[1mm]}
\def\2{\\[2mm]}
\def\Lip{Lipschitz\xspace}
\def\BHP{boundary Harnack principle \xspace}
\def\Note{\nind\textit{Note.}\hskip 2mm}
\def\Remark{\nind\textit{Remark.}\hskip 2mm}
\def\Notation{\nind\textit{Notation.}\hskip 2mm}
\def\Proof{\nind\textit{Proof.}\hskip 2mm}
\newcommand{\tr}[1]{\text{\rm tr}\indx{#1}}
\medskip
\begin{abstract}
We study the generalized boundary value problem for  (E)\; $-\Delta u+|u|^{q-1}u=0$ in a dihedral domain
$\Gw$, when $q>1$ is supercritical. The value of the critical exponent can take only a finite number of values
depending on the geometry of $\Gw$. When $\gm$ is a bounded Borel measure
in a k-wedge, we give necessary and sufficient conditions in order it
be the boundary value of a solution of (E). We also give conditions which ensure that a
boundary compact subset is removable. These conditions are expressed in terms of Bessel capacities
$B_{s,q'}$ in $\BBR^{N-k}$ where $s$ depends on the characteristics of the wedge. This allows us to describe the boundary trace of a positive solution of (E).

\end{abstract}

\maketitle
\tableofcontents

\section{Introduction}
Let $\Gw$ be a bounded  Lipschitz domain in $\BBR^N$, $\gr$  the first eigenfunction of $-\Gd$ in $W^{1,2}_0(\Gw)$ with supremum $1$ and $\gl$  the corresponding eigenvalue, and let $q>1$. A  long-term research on the equation
 \begin{equation}\label{M1}
-\Gd u+|u|^{q-1}u=0 \txt{in}\Gw,
\end {equation}
has been carried out  for more than twenty years by probabilistic and/or analytic methods. Much of the research was focused on three  main problems in domains of class $C^2$:
\vskip 3mm

 \noindent(i) The Dirichlet problem for \eqref{M1} with boundary data given by a finite Borel  measure on $\bdw$.

 \noindent(ii) The characterization of removable singular subsets of $\bdw$ relative to positive solutions of \eqref{M1}.

 \noindent (iii) The characterization of arbitrary positive solutions of \eqref{M1} via an appropriate notion of boundary trace.
\vskip 3mm

Consider the Dirichlet problem
\begin{equation}\label{M1-1}\BA {l}
-\Gd u+|u|^{q-1}u=0 \txt{ in }\Gw,\;
u=\gm \txt{ in }\prt\Gw
\EA\end {equation}
where $\gm\in \GTM(\bdw)$ (= space of finite Borel measures on $\prt\Gw$).  Following \cite{MV7}, a (weak) solution $u:=u_\gm$ of \eqref{M1-1} is a function $u\in L^q_{\gr}(\Gw)$ such that,
\begin{equation}\label{test-1}
\int_{\Gw}\left(-u\Gd\eta+\eta|u|^{q-1}u\right)dx=-\int_{\Gw}\BBK[\gm]\Gd\eta dx,
\end{equation}
for every in $\eta\in X(\Gw)$, where
\begin{equation}\label{test}
X(\Gw)=\left\{\eta\in W^{1,2}_0(\Gw):\gr^{-1}\Gd\eta\in L^\infty(\Gw)\right\}.
\end{equation}
Here $\BBK[\gm]$ is the harmonic function in $\Gw$ with boundary trace $\gm$ and $\rho$ is the first eigenfunction of $-\Gd$ in $\Gw$ normalized so that $max_{_\Gw} \,\rho=1$. We recall that, if $\Gw$ is Lipschitz $\BBK[\gm]\in L^1_{\gr}(\Gw)$; if $\Gw$ is of class $C^2$, $\BBK[\gm]\in L^1(\Gw)$.

A measure $\mu$ is a \emph{$q$-good measure} if \eqref{M1-1} has a solution. The space of $q$-good measures is denoted by $\GTM_q(\bdw)$. It is known that, if $\mu$ is $q$-good, the solution is unique. Furthermore, if $\mu$ satisfies the condition
\begin{equation}\label{admi}
\int_\Gw\BBK[\abs\gm]^q\gr dx<\infty,
\end{equation}
then it is $q$-good. When $\mu$ satisfies this condition we say that it is a \emph{$q$-admissible measure}.

When  $\Gw$ is a domain of class $C^2$, $\BBK[\mu]\in L^q_\rho$ for every $q\in (1,\frac{N+1}{N-1})$ and every $\mu\in \GTM(\bdw)$. Therefore, for $q$ in this range, every measure  in $\GTM(\bdw)$ is $q$-good and there is no removable boundary set (except for the empty set).  Problem (iii), for $q$ in this range, was resolved by Le Gall \cite{LeG} (for $N=q=2$) and Marcus and V\'eron \cite{MV1} (for $1<q<\frac{N+1}{N-1}$, $N\geq 3$).

 The number $q_c=\frac{N+1}{N-1}$ is called the \emph{critical value} for \eqref{M1}. If $q$ is supercritical, i.e. $q\geq q_c$, point singularities are removable. In particular there is no solution of
\eqref{M1-1} when $\mu=\gd_y$ (= a Dirac measure concentrated at a point $y\in \bdw$).

In the supercritical case, problems (i) - (iii), $\Gw$ of class $C^2$, have been resolved in several stages. We say that a compact set $E\sbs \bdw$ is removable relative to equation \eqref{M1} if there exists no positive solution vanishing on $\bdw\sms E$. We say that $E$ is conditionally removable if any solution $u$ of \eqref{M1-1}, with $\mu\in \GTM(\bdw)$, such that $u=0$ on $\bdw\sms E$ must vanish in $\Gw$.

With respect to problem (ii) it was shown that a compact set $E\sbs \bdw$  is removable if and only if
$C_{\frac{2}{q},q'}(E)=0$, $q'=q/(q-1)$. Here $C_{\ga,p}$ denotes the Bessel capacity, with the indicated indexes  on $\bdw$. (see Section 4.2 for an overview of Bessel capacities). This result was obtained by Le Gall \cite{LeG} for $q=2$, Dynkin and Kuznetsov \cite{DK2} for $1<q\leq 2$, Marcus and V\'eron \cite{MV2} for $q>2$. For a unified analytic proof, covering all $q\geq q_c$ see \cite{MV3}.

The above result implies that every $q$-good measure $\mu$ must vanish on sets of $C_{\frac{2}{q},q'}$ capacity zero. On the other hand a result of Baras and Pierre \cite{BP} implies that every positive measure $\mu\in \GTM(\bdw)$
that vanishes on sets of $C_{\frac{2}{q},q'}$ capacity zero is the limit of an increasing sequence of admissible measures and therefore $q$-good. In conclusion: a measure $\mu\in \GTM(\bdw)$ is $q$-good if and only if
it vanishes on sets of $C_{\frac{2}{q},q'}$ capacity zero. This takes care of problem (i).

Problem (iii) has been treated in several papers, with various definitions of a generalized boundary trace
for positive solutions of \eqref{M1}, see \cite{DK3} and \cite{MV6}. Finally a full characterization of positive solutions was obtained by Mselati \cite{Ms} for $q=2$, Dynkin \cite{Dy2} for $1<q<2$ and Marcus \cite{MM} for every $q\geq q_c$. In \cite{Ms, Dy2} the restriction to $q\leq2$ was dictated by their use of probabilistic techniques that do not apply to $q>2$. In \cite{MM} the proof is purely analytic.

\medskip

If $\Gw$ is Lipschitz, $\gx\in \bdw$,  we say that $q_\gx$ is the critical value for \eqref{M1} at $\xi$ if, for $1<q<q_\gx$,
problem \eqref{M1-1} with $\gm=\gd_\xi$ has a solution, but for $q> q_\gx$ no such solution exists.

In contrast to the case of smooth domains,
when $\Gw$ is Lipschitz, $q_\gx$ may vary with the point. For every compact set $F\sbs \bdw$ there exists a number $q(F)>1$ \sth, for $1<q<q(F)$, every measure in $\GTM(\bdw)$ supported in $F$ is $q$-good. Obviously 
$q(F)\leq \min\{q_\gx:\,\gx\in F\}$ but it is not clear if equality holds.


In the special case when  $\Gw$ is a polyhedron,
the function $\gx\to q_\gx$ obtains only a finite number of values
(in fact, it is constant on each open face and each open edge) and, if $q\geq q_\gx$, an isolated singularity at
$\gx$ is removable. Furthermore, the assumption $1<q<\min\{q_\gx:\,\gx\in \bdw\}$ implies that every measure in $\GTM(\bdw)$ is $q$-good. For this and related results see \cite{MV7}.

In the present paper we study problem \eqref{M1-1} when $\Gw$ is a polyhedron and $q$ is supercritical, i.e. $q\geq \min\{q_\gx:\,\gx\in \bdw\}$. Following is a description of the main results.
\medskip

\
\noindent{\bf A.} {\em On the action of Poisson type kernels with fractional
dimension.}

In preparation for the study of supercritical \bvp{}s
 we establish an
harmonic analytic result, extending a well known result on the
action of Poisson kernels on Besov spaces with negative index
(see \cite[1.14.4.]{Tri} and \cite{Qui}). We first quote the
classical result for comparison purposes.

\bprop {pref-repr}Let $1<q<\infty$ and $s>0$. Then, for any bounded
Borel measure $\gm$ in $\BBR^{n-1}$,
\begin{equation}\label{pref-P4}
I(\gm)=\int_
{\BBR_+^{n}}\abs{\BBK_{n}[\gm](y)}^qe^{-y_{1}}y_{1}^{sq-1}dy\approx
\norm{\gm}^q_{B^{-s,q}(\BBR^{n-1})}.
\end{equation}
\es

Here $\BBK_n[\mu]$ denotes the Poisson potential of $\mu$ in
$\BBR^n_+=\BBR_+\ti\BBR^{n-1}$, namely,
\begin{equation}\label{pref-BBKn}
  \BBK_n[\mu](y)=\gg_{n}y_{1}\int_{\BBR^{n-1}}\myfrac{d\gm(z)}{\left(y_1^2+|\gz-z|^2\right)^{n/2}}\forevery y=(y_1,\gz)\in
\BBR^n_{+}
\end{equation}
where $\gg_n$ is a constant depending only on $n$.
\medskip

\noindent\textit{Notation.}\hskip 2mm Let $m$ be a positive integer and let $\nu$ be a real number, $\nu\geq m+1$. Denote,
\begin{equation}\label{pref-Knum}
 \BBK_{\nu,m}[\mu](\tau,\gz):=\int_{\BBR^{m}}
\myfrac{\tau^{\nu-m}d\gm(z)}{(\tau^2+|\gz-z|^2)^{\gn/2}}\forevery
\tau\in (0,\infty),\, \gz\in \BBR^m.\end{equation}
Note that
$$\BBK_n[\mu]=\gg_n\BBK_{n,n-1}[\mu].$$

 \bth{pref-general-nu} Let $m$ and $\nu$ be as above. Then, for every $q>1$ and every
$s\in (0,m/q')$, $q'=q/(q-1)$, there exists a
positive constant $c$ \sth, for every positive measure $\mu\in \GTM(\BBR^m)$ supported in $ B_{R/2}(0)$
for some $R>1$,
\begin{equation}\label{pref-Mns<}\BAL
\rec{c}\norm{\gm}^q_{B^{-s,q}(\BBR^{m})}&\leq
\int_0^R\Big(\int_{|\gz|<R}\abs{\BBK_{\nu,m}[\mu](\tau,\gz)}^q
d\gz\Big) \tau^{ sq-1}\,d\tau\\
& \leq c R^{(s+\nu-m) q+1}\norm{\gm}^q_{B^{-s,q}(\BBR^{m})}.\EAL
\end{equation}

This also holds when $s=m/q'$, provided that the diameter of
$\supp \mu$ is sufficiently small. \es

This is proved in Section 3 (see \rth{general-nu}) using a slightly
different notation.

\medskip

\nind{\bf B.} {\em The critical
value  and the characterization of $q$-good measures in a k-wedge.}

\medskip

The next step towards the study of \bvp{}s in a polyhedron is the
treatment of such problems in  a k-wedge (or k-dihedron) i.e., the
domain defined by the intersection of $k$ hyperplanes in $\BBR^N$,
$1<k<N$. The edge is an $(N-k)$ dimensional space.

We note that if
$k=N$ the 'edge' is a point and the corresponding wedge is a cone with vertex at this point.  If $k=1$ the wedge is a half space. Both of these cases  have been treated in \cite{MV7}.

Let $A$ be a Lipschitz domain in $S^{k-1}$.
If
\begin{equation}\label{pref-dom}
S_A:=\{x\in \BBR^N:\abs x=1,\,x\in A\ti\prod_{j=k}^{N-1}[0,\gp]\}\sbs
S^{N-1}\}
\end{equation}
then
$$D_{A}:=\{x=(r,\gs):r>0,\gs\in S_A\}$$
is a k-wedge in $\BBR^N$ whose `edge' $d_A$ may be identified with
$\BBR^{N-k}$ and its `opening' is $A$.

Let $\gl_{A}$ be the first eigenvalue of $-\Gd_{_{S^{N-1}}}$ in
$W^{1,2}_{0}(S_A)$
 and denote by $\kappa_\pm$ the roots of the equation,
   \begin{equation}\label{kappa1}
\gk^2+(N-2)\gk-\gl_{A}=0.
 \end{equation}
Put
\begin{equation}\label{pref-q-critk}
q_{c}:=\myfrac{\gk_{+}+N}{\gk_{+}+N-2}
\end{equation}
and
\begin{equation}\label {pref-qk}
q^*_c:=1+\myfrac{2-k+\sqrt{(k-2)^2+4\gl_{A}-4(N-k)\gk_{+}}}{\gl_{A}-(N-k)\gk_{+}}.
\end{equation}

Let  $C_{\ga,p}^{N-k}$ denote the Bessel
capacity with the indicated indices in $\BBR^{N-k}$. The next theorem provides a characterization of $q$-good measures supported on $d_A$.

\bth{pref-edge-bvp} (a) If $1<q<q_c$ every measure  in $\GTM(d_A)$ is $q$-good relative to $D_A$. In fact every such measure is $q$-admissible.

\noindent(b) If $q\geq q^*_c$, the only $q$-good measure in $\GTM(d_A)$ is the zero measure.

\noindent(c) If $q_c\leq q<q^*_c$, a measure $\mu\in
\GTM(d_A)$ is  $q$-good
relative to $D_A$ if and only if $\mu$ vanishes on
every Borel set $E\sbs d_A$ \sth $C_{s,q'}^{N-k}(E)=0$,
$s=2-\frac{k+\kappa_+}{q'}$.
\es

 \emph{The characterization of $q$-good measures in a polyhedron} follows as an easy consequence of the above theorem (see \rth{poly-good} below).
\medskip


\nind{\bf C.} {\em Characterization of removable sets.}

\smallskip
Let $\Gw$ be an N-dimensional polyhedron. \rth{pref-edge-bvp}
provides a necessary and sufficient condition for the removability
of a singular set $E$ relative to the family of solutions $u$ \sth
$$\int_\Gw |u|^q\gr\,dx<\infty.$$

 The next
result provides a necessary and sufficient condition for
{\em removability} in the sense that the only non-negative solution $u\in
C(\bar\Gw\sms E)$ which vanishes on $\bar\Gw\sms E$ is the trivial
solution $u=0$.

Let $L$ denote a face or edge or vertex of $\Gw$ and put $k:=\mathrm{codim}\,L$. If $1<k< N$ let $d_L$ denote the linear space spanned by $L$, such that $L$ is an open subset of $d_L$.  Let $Q_L$ denote the k-wedge with boundary $d_L$  such that, for some neighborhood $M$ of $L$, $\Gw\cap M= Q_L\cap M$ and let $A_L$ denote the opening of $Q_L$. If $k=N$, $Q_L$ is a cone with vertex $L$.  Let $q_c(L)$ and $q^*_c(L)$ be defined as in  \eqref{pref-q-critk} and \eqref{pref-qk} for $A=A_L$. Finally let
$$s(L)= 2-\frac{k+\kappa_+}{q'}$$
where $\kappa_\pm$ are the roots of \eqref{kappa1}  for  $A=A_L$.
If $k=N$, $Q_L$ is a cone with vertex $L$. In this case  $q_c(L)=q^*_c(L)=1-\frac{2}{\kappa_-}$. If $k=1$ $q_c(L)=q^*_c(L)=(N+1)/(N-1)$.

\bth{pref-rem} Let $\Gw$ be a polyhedron in $\BBR^N$.
A compact set $E\sbs \bdw$ is removable if and only if, for every  $L$ as above such that $E\cap L\neq \ems$ the following conditions hold.

\noindent If $1\leq k<N$: either
$q_c(L)\leq q<q^*_c(L)$ and  $C^{N-k}_{s(L),q'}(E\cap L)=0$ or $q\geq q^*_c(L)$.

\noindent If $k=N$: $q\geq q_c(L)$.
 \es

\medskip

The present paper is  part of an article, `Boundary trace of positive solutions of semilinear elliptic equations in Lipschitz domains' arXiv:0907.1006 (2009). The first part of this article was published in \cite{MV7}. The second and last part are presented here. The characterization of $q$-good measures, here established in polyhedrons, was recently established in \cite{An_Ma}, for arbitrary Lipschitz domains and a general family of nonlinearities. However  the full removability result, \rth{removable set}, has not been superseded. (In \cite{An_Ma} the authors provided  - in the generality mentioned above -  a characterization of \emph{conditional removability} but not of full removability.) The methods of proof in the two papers are completely different. In the present paper, the characterization of $q$-good measures is based on an extension of a result of \cite{Qui} and \cite[1.14.4.]{Tri} on the action of Poisson kernels on Besov spaces with negative index. In \cite{An_Ma} the proof relies on a relation between elliptic semilinear equations with absorption and linear Schr\"odinger equations.



\section{The Martin kernel and critical values in a k-dimensional dihedron.}

\subsection{The geometric framework}

 An N-dim polyhedra $P$ is a bounded domain bordered by a finite number of hyperplanes.
 Thus the boundary of $P$ is the union of a finite number of sets $\{L_{k,j} : k = 1,\cdots,N, \; j =
1,\cdots, n_k\}$ where
$\{L_{1,j}\}$ is the set of open faces of P,  $\{L_{k,j}\}$ for $k = 2,\cdots,N -1$, is the family of
relatively open N- k-dimensional edges and $\{L_{N,j}\}$ is the family of vertices of $P$.
  An N-k dimensional edge is a relatively open set in the intersection of k hyperplanes; it
  will be described by the characteristic angles of these hyperplanes.
  \medskip

 We recall that the spherical coordinates in $\BBR^N=\{x=(x_{1},...x_{N})\}$ are expressed by
 \begin{equation}\label{N-dim1}\left\{\BA {l}
 x_{1}=r\sin\gth_{N-1}\sin\gth_{N-2}...\sin\gth_{2}\sin\gth_{1}\\
 x_{2}=r\sin\gth_{N-1}\sin\gth_{N-2}...\sin\gth_{2}\cos\gth_{1}\\
 x_{3}=r\sin\gth_{N-1}\sin\gth_{N-2}...\cos\gth_{2}\\
 .\\
 .\\
 .\\
 x_{N-1}=r\sin\gth_{N-1}\cos\gth_{N-2},\\
 x_{N}=r\cos\gth_{N-1}
\EA\right. \end{equation}
where, $r=|x|$,
$\gth_{1}\in [0,2\gp]$ and $\gth_{\ell}\in[0,\gp]$ for $\ell=2,3,...,N-1$. We denote $\gs=(\theta_1,...\theta_{N-1})$. Thus in spherical coordinates $x=(r,\gs)$.
\medskip

We consider an unbounded {\em non-degenerate k-dihedron}, $2\leq k\leq N$ defined as follows. Let
$A$ be given by
 $$A=\begin{cases}(0,\ga_{1})\ti \prod_{j=2}^{k-1}(\ga_{j},\ga'_{j})
 &\text{if $k>2$}\\
(0,\ga_{1})
&\text{if $k=2$}\end{cases}
 $$
 where
 $$0<\ga_{1}< 2\gp,\q 0\leq\ga_{j}<\ga'_{j}<\gp \q j=2,...,k-1.$$
 We denote by $S_A$ the spherical domain
\begin{equation}\label{dom}
S_A=\{x\in \BBR^N:\abs x=1,\,\gs\in
A\ti\prod_{j=k}^{N-1}[0,\gp]\}\sbs S^{N-1}\}
\end{equation}
and by $D_A$ the corresponding  k-dihedron,
$$D_{A}=\{x=(r,\gs):r>0,\gs\in S_A\}.$$
The {\em edge}\/ of $D_{A}$ is the (N-k)-dimensional space
 \begin{equation}\label{edge}
d_{A}=\{x:x_{1}=x_{2}=...=x_{k}=0\}.
\end{equation}
\subsection{On the Martin kernel and critical values in a cone.}\label{ss:cone}
We recall here some elements of local analysis when $\Gw=C_A\cap B_1$, $A$ is a  Lipschitz domain in $S^{N-1}$ and $C_A$ is the cone with vertex $0$ and opening $A$.

Denote by $\gl\indx{A}$ the first eigenvalue and by $\gf\indx{A}$
the first eigenfunction of $-\Gd'$ in $W^{1,2}_{0}(A)$ (normalized by
$\max \gf_{_A}=1$).  Let $\kappa_-$ be the negative root of \eqref{kappa1}
and put
$$\Gf_1(x):=\rec{\gg} |x|^{\kappa_-}\gf_{_A}(x/\abs x)$$
where  $\gg$ is a positive number. Then $\Gf_1$ is a harmonic
function in $C_A$ vanishing on $\prt C_A\sms \{0\}$ . We choose  $\gg=\gg_A$ so that the boundary trace of $\Gf_1$ is $\gd_0$ (=Dirac measure on  with mass $1$ at the origin).
\medskip

\noindent(i) \ If $q\geq 1-\frac{2}{\kappa_-}$,
 there is no solution of \eqref{M1} in $\Gw_S$ with isolated
singularity at $0$. (See \cite{FV}.)
\medskip

\noindent(ii) If $1<q<1-\frac{2}{\kappa_-}$, then for any $k>0$
there exists a unique solution $u:=u_k$ to problem \eqref{M1-1} with
$\gm=k\gd_0$ and

\begin{equation}\label{pref-lim-k}
u_{k}(x)= k\Gf_1(x)(1+o(1))\quad
\txt{as}x\to0.
\end{equation}

The function $u_\infty=\lim_{k\tin} u_k$  is a positive
solution of \eqref{M1} in $\Gw$ which vanishes on $\prt \Gw\sms
\{0\}$ and satisfies
\begin{equation}\label{pref-lim-infi}
u_{\infty}(x)= |x|^{-\frac{2}{q-1}}\gw_{_A}(x/|x|)(1+o(1))\quad
\txt{as}x\to0
\end{equation}
where  $\gw_{_A}$ is the (unique) positive solution of
\begin{equation}\label{pref-NEVP}
-\Gd'\gw-a_{_{N,q}}\gw+\abs\gw^{q-1}\gw=0
\end{equation}
on $S^{N-1}$. Here $\Gd'$ is the Laplace - Beltrami operator and
\begin{equation}\label{pref-NEVP1}
a_{_{N,q}}=\myfrac{2}{q-1}\left(\myfrac{2q}{q-1}-N\right).
\end{equation}

\noindent(iii) If $u\in C(\bar\Gw_{A}\setminus\{0\})$ is a  positive solution of
(\ref{M1}) vanishing on $(\prt C_{A}\cap B_{r_{0}}(0))\sms\{0\}$, then  either $u$ satisfies  \eqref{pref-lim-k} for some $k>0$ or
$u$ satisfies \eqref{pref-lim-infi}. In particular there exists a unique positive solution vanishing on $(\prt C_{A}\cap B_{r_{0}}(0))\sms\{0\}$ with strong singularity at $0$. (For (ii) and (iii) see \cite[Theorem 5.7]{MV7}.)

 \subsection{Separable harmonic functions and the Martin kernel in a k-dihedron, ${2\leq k<N}$.}\hskip 2mm
In the system of spherical coordinates, the Laplacian takes the form
$$\Gd u=\prt_{rr}u+\myfrac{N-1}{r}\prt_{r}u+
\myfrac{1}{r^2}\Gd_{_{S^{N-1}}} u
$$
where the Laplace-Beltrami operator $\Gd_{_{S^{N-1}}}$ is expressed by induction by
 \begin{equation}\label{lapla}\BAL
\Gd_{_{S^{N-1}}} u=&\myfrac{1}{(\sin\gth_{N-1})^{{N-2}}}\myfrac{\prt}{\prt\gth_{N-1}}\left((\sin\gth_{N-1})^{{N-2}}\myfrac{\prt u}{\prt\gth_{N-1}}\right)\\
&+\myfrac{1}{(\sin\gth_{N-1})^{{2}}}\Gd_{_{S^{N-2}}} u.
\EAL \end{equation}
 and
  \begin{equation}\label{lapla1}
\Gd_{_{S^{1}}} u=\prt_{\gth_{1}\gth_{1}}u
 \end{equation}
If we compute the positive harmonic functions in the k-dihedron $D_{A}$ of the form
$$v(x)=v(r,\gs)=r^\gk\gw(\gs) \q\text{in }D_A,\q v=0\q\text{in }\prt D_{A}\setminus\{0\}.
$$
we find that $\gw$ must be a positive eigenfunction corresponding to the first eigenvalue, $\gl_A$, of $-\Gd_{_{S^{N-1}}}$ in $W^{1,2}_{0}(S_A)$,
   \begin{equation}\label{kappa2-0}\left\{\BA {l}
\Gd_{_{S^{N-1}}}\gw+\gl_{A}\gw=0\quad\text{in }S_{A}\\
\phantom{\Gd_{_{S^{N-1}}}\gw+\gl_{A}}\gw=0\quad\text{on }\prt S_A
\EA\right. \end{equation}
and $\gk$ must be a root of the algebraic equation \eqref{kappa1} with $\gl_A$ as above.
 Thus $\kappa=\kappa_\pm$ where
    \begin{equation}\label{kappa2}\BAL
   \gk_{+}&=\myfrac{1}{2}\left(2-N+\sqrt{(N-2)^2+4\gl_{A}}\right)\\
   \gk_{-}&=\myfrac{1}{2}\left(2-N-\sqrt{(N-2)^2+4\gl_{A}}\right).
\EAL\end{equation}
Since
$$
 S^{N-1}=\left\{\gs=(\gs_2\sin\gth_{N-1},\cos\gth_{N-1}):\,\gs_2\in S^{N-2},\;\gth_{N-1}\in (0,\gp)\right\},
$$
we look for a solution $\gw=\gw^{\{1\}}$ of \eqref{kappa2-0} of the form
$$\gw^{\{1\}}(\gs)=(\sin\gth_{N-1})^{\gk_{+}}\gw^{\{2\}}(\gs_2),
\quad\gth_{N-1}\in (0,\gp),\q\gs_2\in S^{N-2}.
$$
Here $S^{N-2}=S^{N-1}\cap\{x_N=0\}$ and we denote
$$S^{\{N-2\}}_A=S_A\cap\{x_{N}=0\},\q D^{\{N-2\}}_{A}:=D_A\cap\{x_{N}=0\}
\subset \BBR^{N-1}.$$
Then (\ref{kappa2})
jointly with relation (\ref{lapla}) implies
\begin{equation}\label{EV2}\left\{\BA {l}
\Gd_{_{S^{N-2}}}\gw^{\{2\}}+(\gl_{A}-\gk_{+})\gw^{\{2\}}=0\quad\text{on }S^{\{N-2\}}_A\\[2mm]
\phantom{\Gd_{_{S^{N-2}}}\gw^{\{2\}}+(\gl_{A}-\gk_{+})}
\gw^{\{2\}}=0\quad\text{on }\prt S^{\{N-2\}}_A.
\EA\right.\end{equation}
Since  we are interested in  $\gw^{\{2\}}$  positive, $\gl_{A}^{\{2\}}:=\gl_{A}-\gk_{+}$ must be
 the first eigenvalue of $-\Gd_{_{S^{N-2}}}$ in $W^{1,2}_{0}(S^{\{N-2\}}_A)$.

Next we look for positive harmonic functions  $\tl u$ in $D^{\{N-2\}}_{A}$ \sth
$$\tilde u(x_1,\ldots,x_{N-1})= r^{\gk'}\gw(\gs_2),\q  \tilde u=0 \text{ on }\prt D^{\{N-2\}}_{A}$$
 The algebraic equation which gives the exponents is
$$(\gk')^2+(N-3)\gk'-\gl_{A}^{\{2\}}=0.
$$
Denote by $\gk'_{+}$  the positive root of this equation.
By the definition of $\gl_{A}^{\{2\}}$,
$$\gk_+^2+(N-3)\gk_+-\gl_{A}^{\{2\}}=\gk_+^2+(N-2)\gk_+-\gl_{A}=0. $$
Therefore $\gk'_{+}=\gk_{+}$.
Accordingly, if $k\geq 3$, we set
$$\gw^{\{2\}}(\gs_2)=(\sin\gth_{N-2})^{\gk_{+}}\gw^{\{3\}}(\gs_{3}),
$$
an find that $\gw^{\{3\}}$ satisfies
\begin{equation}\label{EV3}\left\{\BA {l}
\Gd_{_{S^{N-3}}}\gw^{\{3\}}+(\gl_{A}-2\gk_{+})\gw^{\{3\}}=0\quad\text{in }S^{\{N-3\}}_{A}\\[2mm]
\phantom{\Gd_{_{S^{N-3}}}\gw^{\{3\}}+(\gl_{A}-2\gk_{+})}
\gw^{\{3\}}=0\quad\text{on }\prt S^{\{N-3\}}_{A},
\EA\right.\end{equation}
where
$$S^{\{N-3\}}_{A}=S_A\cap\{x_{N}=x_{N-1}=0\}.$$
Performing this reduction process (N-k) times, we obtain the following results.\smallskip

\noindent(i) If $k>2$ then $\gw=\gw^{N-k}(\gs)$ is given by
\begin{equation}\label{EV4}
\gw(\gs)=(\sin\gth_{N-1}\sin\gth_{N-2}...\sin\gth_{k})^{\gk_{+}}\gw^{\{N-k+1\}}(\gs_{N-k+1})
\end{equation}
where $$\gs_{N-k+1}\in S^{k-1}=S^{N-1}\cap\{x_N=,x_{N-1}=\cdots=x_{k+1}=0\}$$
and $\gw':=\gw^{\{N-k+1\}}$ satisfies
\begin{equation}\label{EV5-0}\left\{\BA {l}
\Gd_{_{S^{k-1}}}\gw'+(\gl_{A}-(N-k)\gk_{+})\gw'=0,\quad\text{in } S^{\{k-1\}}_{A}\\[2mm]
\phantom{\Gd_{_{S^{k-1}}}\gw'+(\gl_{A}-(N-k)\gk_{+})}\gw'=0,\quad\text{on
}\prt S^{\{k-1\}}_{A}, \EA\right.\end{equation}
where
$S^{\{k-1\}}_{A}=S_A\cap\{x_{N}=x_{N-1}=...=x_{k+1}=0\}\approx A$
and $\gl_{A}-(N-k)\gk_{+}$ is the first eigenvalue of the problem.
\medskip

\noindent(ii) If $k=2$ then
\begin{equation}\label{EV5-2}
\gw(\gs)=(\sin\gth_{N-1}\sin\gth_{N-2}...\sin\gth_{2})^{\gk_{+}}
\gw^{\{N-1\}}(\gth_{1})
\end{equation}
where $\gs_{N-1}\in S^{1}\approx \gth_{1}\in (0,2\gp)$, and $\gw^{\{N-1\}}$ satisfies
\begin{equation}\label{EV5-3}\left\{\BA {l}
\Gd_{_{S^{1}}}\gw^{\{N-1\}}+(\gl_{A}-(N-2)\gk_{+})\gw^{\{N-1\}}=0\quad\text{on }S^{\{1\}}_A\\[2mm]
\phantom{\Gd_{_{S^{1}}}\gw^{\{N-1\}}+(\gl_{A}-(N-2)\gk_{+})}\gw^{\{N-1\}}=0\quad\text{on }\prt S^{\{1\}}_A,
\EA\right.\end{equation}
with $\prt S^{\{1\}}_A\approx (0,\ga)$. In this case
\begin{equation}\label{EV6}
\gk_{+}=\myfrac{\gp}{\ga},\q \gw^{\{N-1\}}(\gth_{1})=\sin(\gp\gth_{1}/\ga),
\end{equation}
and, by \eqref{kappa1},
\begin{equation}\label{EV7}
\gl_{A}-(N-2)\gk_{+}=\myfrac{\gp^2}{\ga^2}\Longrightarrow
\gl_{A}=\myfrac{\gp^2}{\ga^2}+(N-2)\myfrac{\gp}{\ga}.
\end{equation}
Observe that $\rec{2}\leq\gk_+$ with equality holding only in the degenerate case $\ga=2\gp$ (which we exclude).
\medskip

In either case,  we find a positive harmonic function $v_{A}$ in
$D_{A}$, vanishing on $\prt D_{A}$, of the form
\begin{equation}\label{vA}
 v_{A}(x)=\abs x^{\gk_{+}}\gw(x/\abs x)
\end{equation}
with $\gw$ as in \eqref{EV4} (for $k>2$) or \eqref{EV6} (for k=2). Furthermore, if  $\Gw$ is a domain in $\RN$
\sth, for some $R>0$,  $\Gw\cap B_R(0)=D_A\cap B_R(0)$ and $w$ is a positive harmonic function in $\Gw$ vanishing on $d_A\cap B_R(0)$ then $w\sim v_A$ in  $\Gw\cap B_{R'}(0)$ for every $R'\in (0,R)$.

Similarly we find a  positive harmonic function in $D_A$ vanishing  on $\prt D_A\sms \{0\}$,
singular at the origin, of the form
$$K'_{A}(x)=\abs x^{\gk_{-}}\gw(x/\abs x).$$
If $\Gw$ is a domain as above and $z$ is a positive harmonic function in $\Gw$ vanishing on $d_A\cap B_R(0)\sms\{0\}$ then $z\sim K'_A$ in  $\Gw\cap B_{R'}(0)\sms\{0\}$ for every $R'\in (0,R)$.

As $K'_A$ is a kernel function of $-\Gd$ at $0$  it follows that $K'_{A}$ is, up to a
multiplicative constant $c_{A}$, the
Martin kernel of $-\Gd$ in $D_{A}$, with singularity at $0$.
  The Martin kernel, with singularity at a
point $z\in d_A$, is given by
\begin{equation}\label{k3}
K_{A}(x,z)=c_{_{A}}\myfrac{(\sin\gth_{N-1}\sin\gth_{N-2}...\sin\gth_{k})^{\gk_{+}}\gw^{\{N-k+1\}}(\gs_{N-k+1})}{|x-z|^{N-2+\gk_{+}}}
\end{equation}
for every $x\in D_A$. From (\ref{N-dim1})
$$\sin\gth_{N-1}\sin\gth_{N-2}...\sin\gth_{k}=|x-z|^{-1}{\sqrt{x_{1}^2+x_{2}^2+...+x_{k}^2}}.
$$
Therefore, if we write $x\in\BBR^N$ in the form $x=(x',x'')$,
$x'=(x_{1},...,x_{k})$, $x''=(x_{k+1},\cdots,x_N)$,
we obtain the formula,
\begin{equation}\label{k4}\BA{l}
K_{A}(x,z)=c_{_{A}}\myfrac{|x'|^{\gk_{+}}\gw^{\{N-k+1\}}(\gs_{N-k+1})}{|x-z|^{(N-2+2\gk_{+})}}\\[4mm]
\phantom{K_{A}(x,z)}
=c_{_{A}}\myfrac{|x'|^{\gk_{+}}\gw^{\{N-k+1\}}(\gs_{N-k+1})}{(|x'|^2+|x''-z|^2)^{(N-2+2\gk_{+})/2}}.
\EA\end{equation}
Therefore, the Poisson potential of a measure
$\gm\in\GTM (d_{A})$ is expressed by
\begin{equation}\label{k5}\BA{l}
\BBK_A[\gm](x)
\quad=c_{_{A}}|x'|^{\gk_{+}}\gw^{\{N-k+1\}}(\gs_{N-k+1})\\[2mm]
\phantom{----------}\ti\myint{\BBR^{N-k}}{}
\myfrac{ d\gm(z)}{(|x'|^2+|x''-z|^2)^{(N-2+2\gk_{+})/2}}.
\EA\end{equation}

\subsection{The admissibility condition}
Consider the boundary value problem
\begin{equation}\label{BVP}\left\{\BA {l}
-\Gd u+\abs u^{q-1}u=0\quad\text {in }D_A\\
\phantom{-\Gd u+\abs u^{q-1}}
u=\gm\in\GTM(\prt D_A).
\EA\right.\end{equation}
Let
\begin{equation}\label{GgR}
\Gg_{R}=\{x=(x',x''):|x'|\leq R,|x''|\leq R\}, \q D_{A,R}:=D_A\cap \Gg_R
\end{equation}
and let $\gr_{R,A}$ denote the first (positive) eigenfunction in
$D_{A,R}:=D_A\cap \Gg_R$. 
In the rest of this section we drop the index $A$ in $K_A$, $\rho_{_{A,R}}$ etc., except for $D_A$, $D_{A,R}$ and $d_A$.

First we observe that a positive Radon measure on $d_A$ is q-good relative to $D_A$  if and only if, for every compact set $F\sbs d_A$, $\mu\,\chr{F}$ is q-good in $D_A$

Now suppose that $\mu$ is compactly supported in $d_A$ and denote its support by $F$. We claim that $\mu$ is q-good in $D_A$ if and only if  it is q-good relative to $D_{A,R}$ for all sufficiently large $R$. Let $R$ be such that $F\sbs B^{N-k}_{R/2}(0)$. Assume that  $\mu$ is q-good in $D_{A,R}$. Let $v_R$ be the solution of \eqref{M1} in $D_{A,R}$ such that $v_R=\mu$ on $d_A\cap \Gg_R$, $v_R=0$ on $\prt D_{A,R}\sms d_A$. Then $v_R$ increases with $R$ and $v=\lim_{R\tin}v_R$ is a solution of \eqref{M1} in $D_A$ with boundary data $\mu$. This proves our claim in one direction; the other direction is obvious.

The condition for $\mu$ to be  q-admissible in $D_{A,R}$ is
\begin{equation}\label{admiss-R}
\int_{D_{A,R}}\BBK^R[|\gm|](x)^q\gr\ssub{R}(x)dx<\infty.
\end{equation}
where $K^R$ is the Martin kernel of $-\Gd$ in $D_{A,R}$. If $R$ is sufficiently large then, in a neighborhood of $F$, $K^R\sim K$ and $\rho^R\sim\rho\sim v_A$.
 Therefore, a sufficient condition for $\mu$ to be q-good in $D_A$ is
\begin{equation}\label{k6}
\int_{\Gg_{R}\cap D_{A}}\BBK[|\gm|](x)|^q\gr(x)dx<\infty \forevery R>0.
\end{equation}
By the first observation in this subsection, it follows that the previous statement remains valid for any positive Radon measure supported on $d_A$.

By \eqref{k3},
\begin{equation}\label{k7}
\BBK[|\gm|](x)\leq \
c_{_{A}}(r')^{\gk_{+}}\int_{\BBR^{N-k}}j(x',x''-z) d|\gm|(z)
\end{equation}
where
\begin{equation}\label{k7bis}
  j(x)=|x|^{-N+2-2\gk_+} \forevery x\in
  \BBR^N.
\end{equation}
 Therefore, using \eqref{vA}, condition \eqref{k6} becomes
\begin{equation}\label{k8}
\int_0^R\int_{|x''|<R} \Big(\int_{\BBR^{N-k}}
j(x',x''-z)d|\gm|(z)\Big)^q (r')^{(q+1)\gk_{+}+k-1}dx''dr'<\infty
\end{equation}
for every
$R>0$.
\subsection{The critical values.}
Relative to the equation
\begin{equation}\label{eqq}
   -\Gd u+|u|^{q-1}u=0
\end{equation}
there exist two thresholds of criticality associated with the edge
$d_A$.

 The first is the
value $q^*_c$ \sth, for $q^*_c\leq q$
 the whole edge $d_{A}$ is removable
 but for $1<q<q^*_c$ there exist non-trivial solutions in $D_A$
 which vanish on $\prt D_A\sms d_A$.
The second $q_c<q^*_c$  corresponds to the removability of  points
on $d_{A}$. For $q\geq q_c$ points on $d_A$ are removable while for
$1<q<q_c$ there exist solutions with isolated point singularities on
$d_A$. In the next two propositions we determine these critical
values.
\bprop{admp0} Assume $q>1$, $1\leq k<N$. Then the condition
\begin{equation}\label {qk}
q<
q^*_c:=1+\myfrac{2-k+\sqrt{(k-2)^2+4\gl_{A}-4(N-k)\gk_{+}}}{\gl_{A}-(N-k)\gk_{+}}
\end{equation}
is necessary and sufficient for the existence of a non-trivial
solution $u$ of (\ref{eqq}) in $D_{A}$ which vanishes on $\prt
D_{A}\sms d_A$. Furthermore, when this condition holds, there exist
non-trivial positive bounded measures $\mu$ on $d_A$ \sth
$\BBK[\gm]\in L^q_{\gr}(\Gg_{R}\cap D_{A})$. \es

\Remark The statement remains true for $k=N$, which is the
case of the cone. In this case $q_c=q^*_c=1-(2/\kappa_-)$ and a straightforward computation yields:
\begin{equation}\label {qkN}
q_c=\frac{N+2+\sqrt{(N-2)^2+4\gl_A}}{N-2+\sqrt{(N-2)^2+4\gl_A}}.
\end{equation}

 \Proof Recall that $\gl_{A}-(N-k)\gk_{+}$ is the first
eigenvalue in $S_A^{\{k-1\}}$ (see \eqref{EV5-0} and the remarks
following it). Let $\gk'_+,\gk'_-$ be the two roots of the equation
$$X^2+(k-2)X-(\gl_{A}-(N-k)\gk_{+})=0,$$
i.e.
$$\gk'_\pm=\rec{2}\big(2-k\pm\sqrt{(k-2)^2+4(\gl_{A}-(N-k)\gk_{+}}\big).$$
Then,  by  \cite[Theorem 5.7]{MV7}, recalled in subsection 2.2, if $1<q<1-(2/{\gk'_-})$
there exists a unique solution of (\ref{eqq}) in the cone
$C_{S_A^{k-1}}$ i.e. the cone with opening $S_A^{k-1}\sbs
S^{k-1}\sbs \BBR^k$
 with trace $a\gd_0$ (where $\gd_0$ denotes the Dirac measure at the vertex of the cone and
 $a>0$). By \eqref{pref-lim-infi} this solution satisfies
 \begin{equation}\label{lim-k'}
u_{a}(x)=a\abs x^{-\ga}\gf(x/\abs x)(1+o(1))\quad\text{as }x\to 0,
\end{equation}
where $\gf$ is the first positive eigenfunction of $-\Gd'$ in
$W^{1,2}_{0}(S_A^{k-1})$ normalized so that $u_1$ possesses trace
$\gd_0$.

The function $u$ given by
$$\tl u_a(x',x'')= u_a(x')\quad\forall (x',x'')\in D_{A}=C_{S_A^{k-1}}\ti \BBR^{N-k},
$$
is a nonzero solution of (\ref{eqq}) in $D_{A}$ which vanishes on
$\prt D_{A}\sms d_A$ and has bounded trace on $d_A$.

A simple calculation shows that $1-(2/{\gk'_-})$ equals $q_c^*$ as
given in \eqref{qk}.

\smallskip

Next, assume that $q\geq q^*_c$ and let $u$ be a solution of
(\ref{eqq}) in $D_{A}$ which vanishes on $\prt D_{A}\sms d_A$.

Given $\ge>0$ let $v_\ge$ be the solution of \eqref{eqq} in
$D_A^{\{N-k-1\}}\setminus \{x'\in \BBR^k:|x'|\leq \ge\}$ \sth

$$v_\ge(x')=\begin{cases} 0, &\txt{if $x'\in \prt D_A^{\{N-k-1\}},\;|x'|>\ge,$}\\
\infty, &\txt{if $|x'|=\ge$.}
\end{cases}$$
Given $R>0$ let $w\ssub{R}$ be the maximal solution in $\{x''\in
\BBR^{N-k}:|x''|<R\}$.

Then the function $u^*$ given by
$$u^*(x',x'')=v_\ge(x')+w\ssub{R}(x'')$$
is a supersolution of \eqref{eqq} in
$D_A\sms\{(x',x''):|x'|>\ge,\;|x''|<R\}$ and it dominates $u$ in
this domain. But $w\ssub{R}(x'')\to 0$ as $R\tin$ and, by \cite{FV},
$v_\ge(x')\to 0$ as $\ge\to 0$. Therefore $u_+=0$ and, by the same
token, $u_-=0$. \qed

\bprop{admpk} Let $A$ be defined as before. Then
\begin{equation}\label{admpk}
  \BBK[\gm]\in L^q_{\gr}(\Gg_{R}\cap D_{A})\forevery \mu\in \GTM(d_A),\forevery R>0
\end{equation}
 if and only if
\begin{equation}\label{q-critk}
1<q<q_{c}:=\myfrac{\gk_{+}+N}{\gk_{+}+N-2}.
\end{equation}
This statement is equivalent to the following:

Condition \eqref{q-critk} is necessary and sufficient in order that
the Dirac measure $\mu=\gd_P$, supported at a point $P\in d_A$, satisfy \eqref{admpk}.
\es \Proof It is sufficient to prove the result relative to the family of measures $\mu$ such that
 $\mu$ is positive,  has compact
support and  $\mu(d_A)=1$. Let $R>1$ be sufficiently large so
that the support of $\mu$ is contained in $\Gg_{R/2}$. The measure
$\mu$ can be approximated (in the sense of weak convergence of
measures) by a \seq $\set{\mu_n}$ of convex combinations of Dirac
measures supported in $d_A\cap \Gg_{R/2}$. For such a \seq
$\BBK[\mu_n]\to \BBK[\mu]$ pointwise and  $\set{\BBK[\mu_n]}$ is
uniformly bounded in $D_A\sms \Gg_{3R/4}$. Therefore it is
sufficient to prove the result when $\gm=\gd_{0}$.  In this case the
admissibility condition \eqref{admi}) is
$$\int_0^R\int_{|x''|<R} j(x)^q
(r')^{(q+1)\gk_{+}+k-1}dx''dr'<\infty,$$ i.e.,
$$\int_0^R\int_0^R
|x|^{q(2-N-2\gk_+)}(r')^{(q+1)\gk_{+}+k-1}(r'')^{N-k-1}dr''dr'<\infty.$$
Substituting $\tau:=r''/r'$ the condition becomes
$$\int_0^R\int_0^{R/r'}
(1+\tau^2)^{\myfrac{q}{2}(2-N-2\gk_+)}(r')^{q(2-N-\gk_+)+\gk_++N-1}\tau^{N-k-1}d\tau\,dr'<\infty.$$
This holds if and only if $q<(\gk_++N)/(\gk_++N-2)$. \qed
\medskip

\noindent\Remark It is interesting to notice that $k$ does not
appear explicitly in \eqref{q-critk}. Furthermore, we  observe that
\begin{equation}\label{q-critk2}
\myfrac{2}{q_{c}-1}\left(\myfrac{2q_{c}}{q_{c}-1}-N\right)=\gl_{A}\Longleftrightarrow \gk_{+}(\gk_{+}+N-2)=\gl_{A},
\end{equation}
which follows from \eqref{kappa2}. This implies that there does not exist a
nontrivial  solution of the nonlinear eigenvalue problem
\begin{equation}\label{nlne1}\BAL
 -\Gd_{_{S^N-1}}\psi-\myfrac{2}{q-1}\left(\myfrac{2q}{q-1}-N\right)\psi+
 |\psi|^{q-1}\psi&=0\q\text{in }S_{_{D_{A}}}\\[2mm]
\psi&=0\q\text{in }\prt S_{_{D_{A}}} \EAL
\end{equation}
which, in turn,  implies that there does not exists a nontrivial solution of
(\ref{eqq}) of the form $u(x)=u(r,\gs)=|x|^{-2/(q-1)}\psi(\gs)$, and
also no solution of this equation in $D_{A}$ which vanishes on $\prt
D_{A}\setminus\{0\}$. This is the classical ansatz for the
removability of isolated singularities in $d_{A}$. 

\section{The harmonic lifting of a Besov space $B^{-s,p}(d_A)$.}
Denote by $W^{\gs,p}(\BBR^\ell)$ ($\gs>0$, $1\leq p\leq\infty$)
the Sobolev spaces over $\BBR^\ell$. In order to use interpolation,
it is useful to introduce the Besov space $B^{\gs,p}(\BBR^{\ell})$
($\gs>0$). If $\gs$ is not an integer then
 \begin{equation}\label{B01}
B^{\gs,p}(\BBR^{\ell})=W^{\gs,p}(\BBR^{\ell}).
\end {equation}
If $\gs$ is an integer the space is defined as follows.
Put
$$\Gd_{x,y}f=f(x+y)+f(x-y)-2f(x).$$
Then
 \begin{equation}\label{B1}
B^{1,p}(\BBR^{\ell})=\left\{f\in L^p(\BBR^{\ell}):
\myfrac{\Gd_{x,y}f}{|y|^{1+\ell/p}}\in L^p(\BBR^{\ell}\ti
\BBR^{\ell})\right\},
\end {equation}
with norm
 \begin{equation}\label{B2}
\norm f_{B^{1,p}}=\norm f_{L^{p}}+ \left(\dint_{\BBR^{\ell}\ti
\BBR^{\ell}}
\myfrac{|\Gd_{x,y}f|^p}{|y|^{\ell+p}}dx\,dy\right)^{1/p},
\end {equation}
(with standard modification if $p=\infty$) and
 \begin{equation}\label{B2'}\BAL
B^{m,p}(\BBR^{\ell})=\Big\{&f\in
W^{m-1,p}(\BBR^{\ell}):\\
&D_x^{\ga}f\in B^{1,p}(\BBR^{\ell})\;\forall \ga\in \BBN^{\ell},\;|\ga|=m-1\Big\}
\EAL
\end{equation}
with norm
 \begin{equation}\label{B2''}\BAL
 \norm f_{B^{m,p}}=\norm f_{W^{m-1,p}}+
\left(\sum_{|\ga|=m-1}\dint_{\BBR^{\ell}\ti \BBR^{\ell}}
\myfrac{|D_x^{\ga}\Gd_{x,y}f|^p}{|y|^{\ell+p}}dx\,dy\right)^{1/p}.
\EAL \end {equation}

We recall that the following inclusions hold (\cite[p 155]{St})
  \begin{equation}\label{B3}\BA {l}
W^{m,p}(\BBR^{\ell})\subset B^{m,p}(\BBR^{\ell})
\quad\text{if }\,p\geq 2\\[2mm]
B^{m,p}(\BBR^{\ell})\subset W^{m,p}(\BBR^{\ell})
\quad\text{if }\,1\leq p\leq 2.
\EA\end {equation}
When $1<p<\infty$, the dual spaces of $W^{s,p}$ and $B^{m,p}$ are respectively denoted by $W^{-s,p'}$
and $B^{-m,p'}$.


The following is the main result of this section.
\bth{main1}Suppose that $q_c<q<q^*_c$ and let $A$ be defined as in
subsection {\bf 2.1}. Then there exist positive constants $c_1,c_2$, depending on $q,N,k,\gk_+$,
 such that for any $R>1$ and any $\gm\in\GTM_+(d_{A})$ with support in $B_{R/2}$:
\begin{equation}\label{L7}\BAL
&c_1\norm{\gm}^q_{B^{-s,q}(\BBR^{N-k})}\\
&\leq \int_{D_{A,R}}\BBK[|\gm|]^q(x)\gr(x)dx \leq
c_2(1+R)^{\gb}\norm{\gm}^q_{B^{-s,q}(\BBR^{N-k})},
\EAL\end{equation}
where $s=2-\frac{\gk_{+}+k}{q'}$, $\gb=(q+1)\gk_++k-1$ and $D_{A,R}=D_{A}\cap\Gg_R$. If $q=q_c$
the estimate
remains valid for measures $\mu$ \sth the diameter of $\supp \mu$ is sufficiently small
(depending on the parameters mentioned before).
 \es
\Remark When $q\geq 2$ the norms in the Besov space may be replaced
by the norms in the corresponding Sobolev spaces.\2
\indent Recall the admissibility condition for a measure
$\gm\in\GTM_{+}(d_{A})$:
$$\int_{D_{A,R}}\BBK[\gm]^q(x)\gr(x)dx <\infty \forevery R>0$$
and the equivalence (see \eqref{k6}--\eqref{k8})
\begin{align}\label{k8'}
&\int_{D_{A,R}}\BBK[\gm]^q(x)\gr(x)dx\approx J^{A,R}(\mu):=\\
&\int_0^R\int_{B''_R} \Big(\int_{\BBR^{N-k}}
\frac{d\gm(z)}{(\tau^2+|x''-z|^2)|)^{(N-2+2\gk_+)/2}}\Big)^q
\tau^{(q+1)\gk_{+}+k-1}dx''d\tau, \notag
\end{align}
where $x=(x',x'')\in
\BBR^k\ti\BBR^{N-k}$,  $\tau=|x'|$ and $B''_R=\set{x''\in\BBR^{N-k}: |x''|<R}$. We denote,
\begin{equation}\label{nu}
   \nu=N-2+2\gk_+.
\end{equation}

If $2\gk_+$ is an integer, it is natural to relate \eqref{k8'} to the Poison potential of $\mu$
in $\BBR_+^n=\BBR_+\ti\BBR_{n-1}$ where $n=N-2+2\gk_+$. We clarify this statement below.

Assuming that
$2\leq n+k-N$, denote
$$y=(y_1,\wtl y,y'')\in \BBR^n,\q \wtl y=(y_2,\cdots,y_{n+k-N}),\q
y''=(y_{n+k-N+1},\cdots,y_n).$$

 The Poisson kernel in $\BBR^n_{+}=\BBR_+\ti\BBR_{n-1}$ is given by
\begin{equation}\label{P1}
P_{n}(y)=\gg_{n}y_{1}|y|^{-n}\quad y_{1}>0,
\end{equation}
for some $\gg_{n}>0$,  and the Poisson potential of a bounded Borel measure $\gm$ with support in
$${\bf d}:=\{y=(0,y'')\in \BBR^n:\,y''\in \BBR_{N-k}\}$$
is
\begin{equation}\label{P2}\BBK_{n}[\gm](y)=
\gg_{n}y_{1}\int_{\BBR^{N-k}}\myfrac{d\gm(z)}{\left(y_1^2+|\wtl y|^2+|y''-z|^2\right)^{n/2}}\forevery y\in \BBR^n_{+}.
\end{equation}
In particular, for $\wtl y=0$,
\begin{equation}\label{P3}\BBK_{n}[\gm](y_{1},0,y'')=\gg_{n}y_{1}\myint{\BBR^{N-k}}{}\myfrac{d\gm(z)}{\left(y_{1}^2+|y''-z|^2\right)^{n/2}}.
\end{equation}
The  integral in \eqref{P3} is precisely the same as the inner
integral in \eqref{k8'}.

In fact, it will be shown that, if we set
\begin{equation}\label{def-n}
n:=\{\nu\}=\inf\set{m\in\BBN: m\geq\nu},
\end{equation}
this approach also works when $2\gk_+$ is not an integer.
We note that, for $n$ given by \eqref{def-n},
\begin{equation}\label{n-N+k}
n-N+k\geq 2,
\end{equation}
with equality only if $k=3$ and $\gk_+\leq 1/2$ or $k=2$ and $\gk_+\in(1/2,1]$.
Indeed, $$n-N+k=k+\{2\gk_+\}-2$$
and (as $\gk_+>0$) $\{2\gk_+\}\geq 1$. If $k=2$ then
$\gk_+>1/2$ and \consy $\{2\gk_+\}\geq 2$. These facts imply our assertion.

We also note that $\gk_+$ is strictly
increasing relative to $\gl_A$ and
\begin{equation}\label{gk+}
\gk_+\begin{cases}=1, &\text{if $D_A=\BBR^N_+$,}\\
<1, &\text{if $D_A\subsetneqq\BBR^N_+$,}\\
>1, &\text{if $D_A\supsetneqq\BBR^N_+$}.\end{cases}
\end{equation}
Finally we observe that  $\gg:=\gl_A-(N-k)\gk_+>0$ (see \eqref{EV5-0}) and, by
\eqref{kappa2} and \eqref{qk}:
\begin{equation}\label{q*c1}
 \gg=\gk_+^2+(k-2)\gk_+,\q q^*_c=1+ \frac{-(k-2)+\sqrt{(k-2)^2+4\gg}}{\gg}.
\end{equation}
Therefore $q^*_c$ is strictly decreasing relative to $\gg$ and \consy also relative to $\gk_+$.

The proof of the theorem is based on the following important result proved in \cite [1.14.4.]{Tri}
\bprop {repr}Let $1<q<\infty$ and $s>0$. Then for any bounded Borel measure $\gm$ in $\BBR^{n-1}$
there holds
\begin{equation}\label{P4}
I(\gm)=\int_
{\BBR_+^{n}}\abs{\BBK_{n}[\gm](y)}^qe^{-y_{1}}y_{1}^{sq-1}dy\approx
\norm{\gm}^q_{B^{-s,q}(\BBR^{n-1})}.
\end{equation}
\es

\smallskip

In the first part of the proof we derive inequalities comparing
$I(\mu)$ and $J^{A,R}(\mu)$. Actually, it is useful to consider  a
slightly more general expression than $I(\mu)$, namely:
\begin{equation}\label{Inu}
I_{\gn,\gs}^{m,j}(\gm):=\int_
{\BBR_+^{m+j}}\abs{\int_{\BBR^{m}}{}\myfrac{y_{1}d\gm(z)}{\left(y_1^2+
|\wtl y|^2+|y''-z|^2\right)^{\gn/2}}}^qe^{-y_{1}}y_{1}^{\gs q-1}dy,
\end{equation}
where $\nu$ is an arbitrary number \sth $\nu>m$, $j\geq 1$ and $\gs>0$. A point $y\in \BBR_+^{m+j}$
is written in the form $y=(y_1,\wtl y, y'')\in \BBR_+\ti \BBR^{j-1}\ti\BBR^m$. We assume that
 $\gm$  is supported in $\BBR^m$. Note that, 
 \begin{equation}\label{Imns}
    I(\mu)=\gg_{n}^qI_{n,s}^{m,j} \txt{where}
  m=N-k,\q j=n-m=n-N+k.
 \end{equation}
Put
\begin{equation}\label{Fnum}
 F_{\nu,m}[\mu](\tau):=\int_{\BBR^{m}}\abs{\int_{\BBR^{m}}
\myfrac{d\gm(z)}{(\tau^2+|y''-z|^2)^{\gn/2}}}^q dy''\forevery
\tau\in [0,\infty).
\end{equation}
With this notation, if $j\geq 2$ then
\begin{equation}\label{InuF}
I_{\gn,\gs}^{m,j}(\gm):=\int_0^\infty\int_
{\BBR^{j-1}} F_{\nu,m}[\mu](\sqrt{y_1^2+|\wtl y|^2}\,)e^{-y_{1}}y_{1}^{(\gs+1) q-1}d\wtl y\,dy_1
\end{equation}
and if $j=1$
\begin{equation}\label{InuF1}
I_{\gn,\gs}^{m,1}(\gm):=\int_0^\infty
 F_{\nu,m}[\mu](y_1)e^{-y_{1}}y_{1}^{(\gs+1) q-1}\,dy_1
\end{equation}


\blemma {repl}
Assume that $m<\gn$, $0<\gs$, $2\leq j$ and
$1<q<\infty$.
Then there exists a positive constant $c$, depending on $m,j,\nu,\gs,q$, such
that, for every bounded Borel measure $\gm$ with support in $\BBR^m$:

\begin{equation}\label{Ck2}
\rec{c}\int_0^\infty\,F_{\nu,m}[\mu](\tau)h_{\gs,j}(\tau)d\tau\leq
\,I_{\gn,\gs}^{m,j}(\gm)\leq
c\int_0^\infty\,F_{\nu,m}[\mu](\tau)h_{\gs,j}(\tau)d\tau,
\end{equation}
where $F_{\nu,m}$ is given by \eqref{Fnum} and, for every $\tau>0$,
\begin{equation}\label{hsj}
  h_{\gs,j}(\tau)=\begin{cases}\myfrac{\tau^{(\gs+1)q+j-2}}{(1+\tau)^{(\gs+1)q}}, &\text{if $j\geq2$,}\\[4mm]
  e^{-\tau}\tau^{(\gs+1)q-1},  &\text{if $j=1$.}\end{cases}
\end{equation}

\es
\Proof There is nothing to prove in the case $j=1$. Therefore we assume that $j\geq2$.

We use the notation $y=(y_{1},\wtl
y,y'')\in\BBR\ti\BBR^{j-1}\ti\BBR^{m}$. 
The integrand in
 \eqref{InuF} depends only on $y_1$ and $\gr:=|\wtl y|$. Therefore,
 $I_{\gn,\gs}^{m,j}$ can be written in the form
$$\BAL
I_{\gn,\gs}^{m,j}(\gm)
=c_{m,j}\int_0^\infty\int_0^\infty
F_{\nu,m}[\mu](\sqrt{y_1^2+\gr^2})e^{-y_1}y_1^{(\gs
+1)q-1}\,dy_1\gr^{j-2}d\gr.\EAL
$$
We substitute
 $y_1=(\tau^2-\gr^2)^{1/2}$, then  change the order of integration
 and finally  substitute $\gr=r\tau$. This yields,
$$\BAL &c_{m,j}^{-1}I_{\gn,\gs}^{m,j}(\gm)\\
&=\int_0^\infty\int_\gr^\infty
F_{\nu,m}[\mu](\tau)\gr^{j-2}e^{-\sqrt{\tau^2-\gr^2}}(\tau^2-\gr^2)^{(\gs +1)q/2-1}\tau\,d\tau\,d\gr\\
&=\int_0^\infty\int_0^\tau\,F_{\nu,m}[\mu](\tau)\gr^{j-2}e^{-\sqrt{\tau^2-\gr^2}}(\tau^2-\gr^2)^{(\gs
+1)q/2-1}\tau\,d\gr\,d\tau\\
&=\int_0^\infty\int_0^1\,F_{\nu,m}[\mu](\tau)\tau^{j-2+(\gs
+1)q}e^{-\tau\sqrt{1-r^2}}f(r) dr\,d\tau,
 \EAL$$
 where
 $$f(r)= r^{j-2}(1-r^2)^{(\gs +1)q/2-1}.$$
We denote
$$I_\gs^j(\tau)=\int_0^1\,e^{-\tau\sqrt{1-r^2}}f(r) dr,$$
so that
\begin{equation}\label{Isn}
 I_{\gn,\gs}^{m,j}(\gm)=c_{m,j}\int_0^\infty F_{\nu,m}[\mu](\tau)\tau^{j-2+(\gs
+1)q}I_\gs^j(\tau)d\tau.
\end{equation}

To complete the proof we estimate $I^j_\gs$. Since $j\geq2$,
$f\in L^1(0,1)$  and $I^j_\gs$ is continuous in
$[0,\infty)$ and positive everywhere. Hence, for every $\ga>0$,
there exists a positive constant $c_\ga=c_\ga(\gs)$ \sth
\begin{equation}\label{I*1}
\rec{c_\ga}\leq I_\gs^j\leq c_\ga \txt{in} [0,\ga).
\end{equation}

Next we estimate $I^j_\gs$ for large $\tau$. Since $j\geq2$,
$$I^j_\gs\leq 2^{(\gs +1)q/2-1}\int_0^1(1-r)^{(\gs +1)q/2-1}e^{-\tau\sqrt{1-r}}dr.$$
Substituting $r=1-t^2$ we obtain,
\begin{equation}\label{I*2'}
 I^j_\gs\leq 2^{(\gs +1)q/2}\int_0^1 t^{(\gs +1)q-1}e^{-t\tau}dt=c(\gs,q)\tau^{-(\gs +1)q}.
\end{equation}

On the other hand, if $\tau\geq 2$,
\begin{equation}\label{I*3}\BAL
I^j_\gs(\tau)&= \int_0^1(1-t^2)^{(j-3)/2}t^{(\gs+1)q-1} e^{-\tau t}dt\\
 &=\tau^{-(\gs+1)q} \int_0^\tau
 (1-(s/\tau)^2)^{(j-3)/2}s^{(\gs+1)q-1}e^{-s}ds\\
 &\geq \tau^{-(\gs+1)q}2^{-(j-3)}
 \int_0^1s^{(\gs+1)q-1}e^{-s}ds.
 \EAL\end{equation}
Combining \eqref{Isn} with \eqref{I*1}--\eqref{I*3} we obtain
\eqref{Ck2}. \qed

\par Next we derive an estimate in which integration over $\BBR^n_+=\BBR_+^j\ti\BBR^m$
is replaced by integration over a bounded domain, for measures supported
in  a fixed bounded subset of $\BBR^m$.

Let $B_R^j(0)$ and $B_R^m(0)$ denote the balls of radius $R$
centered at the origin, in $\BBR^j$ and $\BBR^m$ respectively. Denote
\begin{equation}\label{FnumR}
  F^R_{\nu,m}[\mu](\tau)=\int_{B^m_R}\abs{\int_{\BBR^{m}}
\myfrac{d\gm(z)}{(\tau^2+|y''-z|^2)^{\gn/2}}}^q dy''\forevery
\tau\in [0,\infty)
\end{equation}
and, if $j\geq 2$,
\begin{equation}\label{InuR}\BAL
I_{\gn,\gs}^{m,j}(\gm;R)=\int_{B^j_R\cap\{0<y_1\}}
F^R_{\nu,m}[\mu](\sqrt{y_1^2+|\wtl y|^2}\,)e^{-y_{1}}y_{1}^{\gs q-1}d\wtl y\,dy_1.
\EAL\end{equation}
where $(y_1,\wtl y)\in \BBR\ti\BBR^{j-1}$. If $j=1$ we denote,
\begin{equation}\label{InuR1}\BAL
I_{\gn,\gs}^{m,1}(\gm;R)=\int_0^R
F^R_{\nu,m}[\mu](y_1)e^{-y_{1}}y_{1}^{\gs q-1}\,dy_1.
\EAL\end{equation}

Similarly to \rlemma{repl} we obtain,
\blemma{replR} If $j\geq1$,
there exists a positive constant $c$ such
that, for any bounded Borel measure $\gm$ with support in $\BBR^m\cap B_R$

\begin{equation}\label{Ck2R}
c^{-1}\int_0^R\,F^R_{\nu,m}[\mu](\tau)h_{\gs,j}(\tau)d\tau\leq
\,I^{m,j}_{\gn,\gs}(\gm;R)\leq
c\int_0^R\,F^R_{\nu,m}[\mu](\tau)h_{\gs,j}(\tau)d\tau
\end{equation}
with $h_{\gs,j}$ as in \eqref{hsj}.
\es

\nind\Proof In the case $j=1$ there is nothing to prove .
Therefore we assume that $j\geq 2$.

From \eqref{InuR} we obtain,
$$\BAL
I^{m,j}_{\gn,\gs}(\gm;R)=c_{m,j}\int_0^R\int_0^{\sqrt{R^2-\gr^2}}
F^R_{\nu,m}[\mu](\sqrt{y_1^2+\gr^2})e^{-y_1}y_1^{(\gs
+1)q-1}dy_1\gr^{j-2}d\gr.\EAL
$$
Substituting $y_1=(\tau^2-\gr^2)^{1/2}$, then  changing the order of integration
 and finally  substituting $\gr=r\tau$ we obtain,

$$c_{m,j}^{-1}I^{m,j}_{\gn,\gs}(\gm;R)=\int_0^R\int_0^1\,F^R_{\nu,\mu}[\mu](\tau)\tau^{j-2+(\gs
+1)q}e^{-\tau\sqrt{1-r^2}}f(r) dr\,d\tau.$$
 where
 $$f(r)= r^{j-2}(1-r^2)^{(\gs +1)q/2-1}.$$
The remaining part of the proof is the same as for \rlemma{repl}.
\qed
\blemma{reduc}Let $1<q$, $0<\gs$ and assume that $m<\nu q$ and
$0\leq j-1<\nu$.
Then there exists a positive constant $\bar c$, depending on $j,m,q,\gs,\nu$, such
that, for every $R\geq1$ and every bounded Borel measure $\gm $ with support in $B_{R/2}(0)\cap
\BBR^m$,
\begin{equation}\label{R1}\BAL
\abs{\int_0^\infty\,F_{\nu,m}[\mu](\tau)h_{\gs,j}(\tau)d\tau-
\int_0^R\,F_{\nu,m}^R[\mu](\tau)h_{\gs,j}(\tau)d\tau}&\\[2mm]
\leq \bar c R^{(\gs +1-\nu)q+m+j-1}\norm{\gm}_{\GTM}^q&
\EAL\end{equation}
 with $h_{\gs,j}$ as in \eqref{hsj}. \es
\Proof
We estimate,
\begin{equation}\label{J1+J2}
\BAL
&\abs{\int_0^\infty\,F_{\nu,m}[\mu](\tau)h_{\gs,j}(\tau)d\tau-\int_0^R\,F_{\nu,m}^R[\mu](\tau)h_{\gs,j}(\tau)d\tau}\leq\\
&\int_R^\infty \abs{F_{\nu,m}[\mu]}(\tau)h_{\gs,j}(\tau)d\tau+
\int_0^R\,\abs{F_{\nu,m}[\mu]-F_{\nu,m}^R[\mu]}(\tau)h_{\gs,j}(\tau)d\tau.
\EAL
\end{equation}
For every $\tau>0$,
\begin{equation}\label{Fnm<}
\abs{F_{\nu,m}[\mu]}(\tau)\leq \tau^{-\nu q}\norm{\mu}_{\GTM}^q.
\end{equation}
Since $j-1<\nu q$, it follows that
\begin{equation}\label{J1-temp}
\BAL \int_R^\infty \abs{F_{\nu,m}[\mu]}(\tau)h_{\gs,j}(\tau)d\tau
&\leq \norm{\mu}_{\GTM}^q\int_R^\infty \tau^{-\nu q}h_{\gs,j}(\tau)d\tau\\
&\leq c(\gs,q)\norm{\mu}_{\GTM}^q \int_R^\infty \myfrac{\tau^{(\gs+1)q+j-2-\nu q}}{(1+\tau)^{(\gs+1)q}}d\tau\\
&\leq \frac{c(\gs,q)}{\nu q-j+1}
\norm{\mu}_{\GTM}^q R^{j-1-\nu q}.\EAL
\end{equation}
Since, by assumption, $\supp\mu\sbs B_{R/2}$, we have
\begin{equation}\label{J2-temp}
\BAL \int_0^R\,&\abs{F_{\nu,m}[\mu]-F_{\nu,m}^R[\mu]}(\tau)h_{\gs,j}(\tau)d\tau\\
&\leq\int_0^R\,\int_{|y''|>R}\abs{\int_{\BBR^{m}}
\myfrac{d\gm(z)}{(\tau^2+|y''-z|^2)^{\gn/2}}}^q dy''h_{\gs,j}(\tau)d\tau\\
&\leq \norm{\mu}_{\GTM}^q \int_0^R\int_{|\gz|>R/2}(|\tau^2+|\gz|^2)^{-\nu q/2}\,d\gz\,h_{\gs,j}d\tau\\
&\leq c(m,q)\norm{\mu}_{\GTM}^q \int_0^R\int_{R/2}^\infty(\tau^2+\gr^2)^{-\nu q/2}\gr^{m-1}\,d\gr\,h_{\gs,j}d\tau\\
&\leq c(m,q)\norm{\mu}_{\GTM}^q \int_0^R\,\tau^{m-\nu q}\int_{R/2\tau}^\infty(1+\eta^2)^{-\nu q/2}\eta^{m-1}\,d\eta\,h_{\gs,j}d\tau\\
&\leq \frac{c(m,q)}{\nu q-m}\norm{\mu}_{\GTM}^q R^{m-\nu q}\int_0^R\tau^{(\gs+1)q+j-2}\,d\tau\\
&\leq \frac{c(m,q)}{(\nu q-m)((\gs+1)q+j-1)}\norm{\mu}_{\GTM}^q R^{(\gs+1)q+j-1+m-\nu q}.\EAL
\end{equation}
Combining \eqref{J1+J2}--\eqref{J2-temp} we obtain
 \eqref{R1}.
\qed
\bcor{reduc}
For every $R>0$ put
\begin{equation}\label{Jns}
  J_{\nu,\gs}^{m,j}(\mu;R):=\int_0^RF^R_{\nu,m}[\mu](\tau) \tau^{(\gs+1)q+j-2}d\tau.
\end{equation}
Then
\begin{equation}\label{Est1-Jns}\BAL
\rec{c}I_{\gn,\gs}^{m,j}(\gm)-\bar c R^{\gb}\norm{\gm}_{\GTM}^q\leq
J_{\nu,\gs}^{m,j}(\mu;R)
\leq c R^{(\gs+1)q}I_{\gn,\gs}^{m,j}(\gm),\\
\gb=(\gs +1-\nu)q+j+m-1, \EAL
\end{equation}
for every $R>1$ and  every bounded Borel measure $\gm $ with support in $B^m_{R/2}(0):=B_{R/2}(0)\cap
\BBR^m$.
\es
\Proof This is an immediate consequence of \rlemma{reduc} and \rlemma{repl}.
\qed
\blemma{Est2-Jns} Let $m,j$ be positive integers \sth $j\geq 1$ and
let $1<q$, $0<\gs$. Put $n:=m+j$.

Then there exist positive constants $c,\bar c$, depending on $j,m,q,\gs$, such
that, for every $R>1$ and every measure $\gm\in \GTM_+(B^m_{R/2}(0))$,
\begin{equation}\label{Est2-Jns}\BAL
&\rec{c}\norm{\gm}^q_{B^{-\gs,q}(\BBR^{n-1})}-\bar c R^{q\left(\gs-\frac{n-1}{q'}\right)}\norm{\gm}_{\GTM}^q\leq
J_{n,\gs}^{m,j}(\mu;R)\\
&\leq c R^{(\gs+1)q}\norm{\gm}^q_{B^{-\gs,q}(\BBR^{n-1})}.\EAL
\end{equation}

If  $\gs<\frac{n-1}{q'}$,  there exists $R_0>1$ \sth, for all $R>R_0$
\begin{equation}\label{Est3-Jns}\BAL
\rec{2c}\norm{\gm}^q_{B^{-\gs,q}(\BBR^{n-1})}\leq
J_{n,\gs}^{m,j}(\mu;R).
\EAL
\end{equation}
If $\gs=\frac{n-1}{q'}$ then,  there exists $a>0$ \sth the inequality remains valid
for measures $\mu$ \sth $\mathrm{diam}(\supp\mu)\leq a$.

If, in addition, $\frac{j-1}{q'}<\gs$ then
\begin{equation}\label{Est4-Jns}\BAL
\rec{2c}\norm{\gm}^q_{B^{-s,q}(\BBR^{m})}\leq
J_{n,\gs}^{m,j}(\mu;R)
\leq c R^{(\gs+1)q}\norm{\gm}^q_{B^{-s,q}(\BBR^{m})},\EAL
\end{equation}
where $s:=\gs-\frac{j-1}{q'}$.
\es
\Remark Assume that $\mu\geq0$. Then:\\
 (i) If $\mu\in B^{-\gs,q}(\BBR^{n-1})$ and $\frac{j-1}{q'}\geq \gs$ then $\mu(\BBR^m)=0$.\\
(ii) If $\mu\in B^{-s,q}(\BBR^m)$ and $\gs> (n-1)/q'$ then $s> m/q'$ and therefore
$B^{s,q'}(\BBR^m)$ can be embedded in $C(\BBR^m)$.\2
\Proof Inequality \eqref{Est2-Jns} follows from \eqref{Est1-Jns} and \rprop{repr} (see also \eqref{Imns}).

For positive measures $\mu$,
$$\norm{\gm}_{\GTM}=\mu(\BBR^{n-1})\leq \norm{\gm}^q_{B^{-\gs,q}(\BBR^{n-1})}.$$
Therefore, if $\gs<\frac{n-1}{q'}$,  \eqref{Est2-Jns}   implies that
 there exists $R_0>1$ \sth
\eqref{Est3-Jns} holds for all $R>R_0$.

If $\gs=\frac{n-1}{q'}$ \eqref{Est2-Jns}   implies that
$$\rec{c}\norm{\gm}^q_{B^{-\gs,q}(\BBR^{n-1})}-\bar c \norm{\gm}_{\GTM}^q\leq
J_{n,\gs}^{m,j}(\mu;R).$$
But if $\mu$ is a positive bounded measure \sth $\mathrm{diam}(\supp\mu)\leq a$ then
$$\norm{\gm}_{\GTM}/\norm{\gm}^q_{B^{-\gs,q}(\BBR^{n-1})}\to 0 \txt{as $a\to 0$.}$$
The last inequality follows from the
imbedding theorem for Besov spaces according to which
there exists a
continuous trace operator $T:B^{\gs,q'}(\BBR^{n-1})\mapsto
B^{s,q'}(\BBR^{m})$ and a continuous lifting
$T':B^{s,q'}(\BBR^{m})\mapsto B^{\gs,q'}(\BBR^{n-1})$ where $s=\gs-\frac{n-m-1}{q'}$.
\qed

If $\nu\in\BBN$ and $\gs=s+\frac{\nu-m-1}{q'}$,
\begin{equation*}\BAL
 J_{\nu,\gs}^{m,\nu-m}(\mu;R)&= \int_0^RF^R_{\nu,m}[\mu](\tau) \tau^{(\gs+1)q+\nu-m-2}\,d\tau\\
& = \int_0^RF^R_{\nu,m}[\mu](\tau) \tau^{(s+\nu-m)q-1}\,d\tau.
\EAL\end{equation*}
However, if $\mu$ is positive, the expression
\begin{equation}\label{Mns}
M_{\nu,s}^{m}(\mu;R):=\int_0^RF^R_{\nu,m}[\mu](\tau) \tau^{(s+\nu-m) q-1}\,d\tau,
\end{equation}
is meaningful for any real $\nu>m$ and $s>0$. Furthermore, as shown  below, the results stated in \rlemma{Est2-Jns}
can be extended to this general case.
\bth{general-nu} Let $1<q$, $\nu\in \BBR$ and $m$ a positive
integer. Assume that $1\leq\nu-m$ and $0<s<m/q'$. Then there
exists a positive constant $c$ \sth, for every bounded positive
measure $\mu$ supported in $\BBR^m\cap B_{R/2}(0)$, $R>1$,
\begin{equation}\label{Est-Mns<}\BAL
\rec{c}\norm{\gm}^q_{B^{-s,q}(\BBR^{m})}\leq M_{\nu,s}^{m}(\mu;R)
\leq c R^{(s+\nu-m) q+1}\norm{\gm}^q_{B^{-s,q}(\BBR^{m})}.\EAL
\end{equation}
This also holds when $s=m/q'$, provided that the diameter of $\supp
\mu$ is sufficiently small. \es
\Proof If $\nu$ is an integer and $j:=\nu-m$ then this statement is part of \rlemma{Est2-Jns}.
Indeed the condition $s>0$ means that $\gs=s+\frac{j-1}{q'}>\frac{j-1}{q'}$ and the condition
$s< m/q' $ means that $\gs<\frac{n-1}{q'}$.

Therefore we assume that
$\nu\nin\BBN$. Let $n:=\{\nu\}$ and $\gth:=n-\nu$ so that $0<\gth<1$. 
Our assumptions imply that $1\leq n-m-1$ because (as $\nu$ is not an integer) $\nu-m>1$ and \consy
$n-m\geq 2$.

 If $a,b$ are positive numbers, put
\begin{equation*}\label{abt}
   A_\nu:=\frac{a^{(s+\nu-m)q-1}}{(a^2+b^2)^{\nu q/2}}.
\end{equation*}
Obviously $A_\nu$ decreases as $\nu$ increases. Therefore, $A_n\leq A_\nu\leq A_{n-1}$ which in turn implies,
$$M_{n,s}^{m}\leq M_{\nu,s}^{m}\leq M_{n-1,s}^{m}.$$
By \rlemma{Est2-Jns}, the assertions of the theorem are valid in the case that $\nu=n$ or $\nu=n-1$.
Therefore the previous
inequality implies that the assertions hold for any real $\nu$ subject to the conditions imposed.
\qed

By \eqref{k8'},
$$J^{A,R}=\int_0^R F_{\nu,m}^R(\tau)\tau^{(q+1)\gk_{+}+k-1}d\tau,$$
where $m=N-k$ and $\nu=N-2+2\gk_+$. 
\Consy, by \eqref{Jns},
$$J^{A,R}=M_{\nu,s}^m$$
where $s$ is determined by,
$$(s+\nu-m) q-1=(q+1)\gk_{+}+k-1,\q k=\nu-m+2-2\gk_+.$$
It follows that
$$sq=-(k-2+2\gk_+)q+(q+1)\gk_{+}+k=k(1-q)+2q-\gk_+(q-1)$$
and therefore
$$s=2-\frac{k+\gk_+}{q'}.$$

\medskip




\medskip

\noindent{\it Proof of \rth{main1}.}

Put
\begin{equation}\label{nusm}
 \nu:=N-2+2\gk_+,\q  s:=2-\frac{\gk_{+}+k}{q'},\q m:=N-k.    
\end{equation}

Recall that in the case $k=2$ we have $\gk_+>1/2$. Therefore
\begin{equation}\label{nu>m-1}
  \nu-m-1=k-3+2\gk_+>0.
\end{equation}
Furthermore,
$$(s+\nu-m) q-1=(q+1)\gk_{+}+k-1,\q k=\nu-m+2-2\gk_.$$
Thus
$$J^{A,R}=\int_0^R F_{\nu,m}^R(\tau)\tau^{(q+1)\gk_{+}+k-1}d\tau=M_{\nu,s}^m.$$

Next we show that $ 0<s\leq m/q'$.
More precisely we prove
\begin{equation}\label{0<s<m/q'}
 0<s\leq m/q' \;\iff\;q_c\leq q<q^*_c.
\end{equation}

Let $\mu$ be a bounded non-negative Borel measure in $B^{-s,q}(\BBR^m)$. If $s\leq 0$,
$B^{-s,q}(\BBR^m)\sbs L^q(\BBR^m)$. Therefore, in this case, every bounded Borel measure
on $\BBR^m$ is admissible i.e. satisfies \eqref{admpk}. \Consy, by \rprop{admpk}, $q<q_c$.
As we assume $q\geq q_c$ it follows that $s>0$.

If, $s>0$ and $sq'-m\geq0$ then $C_{s,q'}(K)=0$
for every compact subset of $\BBR^m$ and \consy $\mu(K)=0$ for any such set.
Conversely, if $sq'-m<0$ then there exist non-trivial positive bounded measures in $B^{-s,q}(\BBR^m)$.
Therefore, by \rprop{admp0},  $sq'<m$ if and only if $q<q^*_c$.

In conclusion, $0<s\leq m/q'$ and $\nu-m\geq 1$; therefore \rth{main1} is a consequence of
\rth{general-nu}.
\qed

\noindent\Remark Note that the critical exponent for the imbedding
of $B^{2-\frac{\gk_{+}+k}{q'},q'}(\BBR^{N-k})$ into $C(\BBR^{N-k})$
is again
$$ q=q_{c}=\myfrac{N+\gk_{+}}{N+\gk_{+}-2}.
$$


\section{Supercritical equations in a polyhedral domain}

In this section $q$ is a real number larger than $1$ and $P$ an
N-dim polyhedral domain as described in subsection 6.1. Denote by
$\set{L_{k,j}:k=1,\dots,N,\;j=1,\dots, n_k}$ the family of faces,
edges and vertices of $P$. In this notation, $L_{1,j}$ denotes one
of the open faces of $P$; for $k=2,\dots,N-1$, $L_{k,j}$ denotes a
relatively open $N-k$-dimensional edge and $L_{N,j}$ denotes a
vertex. For $1\leq k<N$, the $(N-k)$ dimensional space which
contains $L_{k,j}$ is denoted by $\BBR^{N-k}_j$. If $1<k<N$, the
cylinder of radius $r$ around the axis $\BBR^{N-k}_j$ will be
denoted by $\Gg^\infty_{k,j,r}$ and  the subset $A_{k,j}$ of
$S^{k-1}$ is defined by
$$\lim_{r\to0}\rec{r}(\prt\Gg^\infty_{k,j,r}\cap P)=L_{k,j}\ti A_{k,j}.$$
$A_{k,j}$ is the 'opening' of $P$ at the edge $L_{k,j}$.
For $k=N$ we replace in this definition the cylinder $\Gg^\infty_{N,j,r}$ by the ball $B_r(L_{N,j})$.
For $1<k\leq N$ and $A=A_{k,j}$ we use  $d_{A}$ as an alternative notation for $\BBR^{N-k}_j$ and denote by $D_A$ the k-dihedron  with edge $d_A$ and opening $A$ as in subsection 6.1 (with  $S_A$  defined as in \eqref{dom}). For $k=1$,
$D_A$ stands for the half space $\BBR^{N-1}_j\ti (0,\infty)$.

\subsection{Definitions and auxiliary results.}  Let $\Gw$ be a bounded Lipschitz domain.
We say that  $\set{\Gw_n}$ is a {\em
\Lip exhaustion} of $\Gw$ if, for every $n$, $\Gw_n$
is \Lip and
\begin{equation}\label{exhaustion}
  \Gw_n\sbs \bar \Gw_n\sbs \Gw_{n+1},\q \Gw=\cup\Gw_n,\q
  \BBH_{N-1}(\bdw_n)\to \BBH_{N-1}(\bdw).
\end{equation}

If $\gw_n$ (respectively $\gw$) is the harmonic measure in $\Gw_n$ (respectively $\Gw$)  relative to $x_0\in \Gw_1$, then, for every $Z\in C(\bar \Gw)$,
\begin{equation}\label{hm-w-conv}
   \lim_{n\to\infty}\int_{\bdw_n}Z\,d\gw_n=\int_{\bdw}Z\,d\gw.
\end{equation}
\cite[Lemma 2.1]{MV7}. Furthermore, if $\gm$ is a bounded Borel measure on $\bdw$ and $v:=\BBK^\Gw[\mu]$, there holds
\begin{equation}\label{tr-w-conv}
   \lim_{n\to\infty}\int_{\bdw_n}Zv\,d\gw_n=\int_{\bdw}Z\,d\mu,
\end{equation}
\cite[Lemma 2.2]{MV7}. If $v$ is a positive solution and \eqref{tr-w-conv} holds we say that $\mu$ is the \textit{boundary trace} of $v$.

 The following estimates are proved in \cite[Lemma 2.3]{MV7}



\bprop{K(mu)} Let $\mu$ be bounded Borel measures on $\bdw$).
Then $\BBK[\mu]\in L^1_\gr(\Gw)$ and there exists a constant $C=C(\Gw)$ \sth
\begin{equation}\label{K(mu)<}
   \norm{\BBK[\mu]}_{L^1_\gr(\Gw)}\leq C\norm{\mu}_{\GTM(\bdw)}.
\end{equation}
In particular if $h\in L^1(\bdw;\gw)$ then
\begin{equation}\label{K(h)<}
   \norm{\BBP [h]}_{L^1_\gr(\Gw)}\leq C\norm{h}_{L^1(\bdw;\gw)}.
\end{equation}
\es

The nest result will be used in deriving estimates in a k-dimensional dihedron when the boundary data is concentrated on the edge.

\bprop{IBP} We denote by $G^{\Gw_n}$ (respectively $G^{\Gw}$) the Green function in $\Gw_n$ (respectively $\Gw$). Let $v$ be a positive harmonic function in $\Gw$ with
boundary trace $\gm$. Let  $Z\in C^2(\overline\Gw)$  and let $\tilde
G\in C^\infty(\Gw)$ be a function that coincides with $x\mapsto
G(x,x_{0})$ in $Q\cap \Gw$ for some neighborhood $Q$ of $\prt\Gw$
and some fixed $x_{0}\in\Gw$. In addition assume that there exists a
constant $c>0$ \sth
\begin{equation}\label{ibp0}
|\nabla Z\cdot\nabla \tl G|\leq c\gr.
\end{equation}
Under these assumptions, if $\gz:=Z\tilde G$ then
\begin{equation}\label{ibp1}
-\myint{\Gw}{}v\Gd\gz\, dx=\myint{\prt\Gw}{}Z d \gm.
\end{equation}
\es
 \Proof Let $\set{\Gw_n}$ be a $C^1$ exhaustion of $\Gw$. We assume
that $\bdw_n\sbs Q$ for all $n$ and $x_0\in \Gw_1$. Let $\tl G_n(x)$
be a function in $C^1(\Gw_n)$ \sth $\tl G_n$ coincides with
$G^{\Gw_n}(\cdot, x_0)$ in $Q\cap\Gw_n$, $\tl G_n(\cdot,x_0)\to \tl
G(\cdot,x_0)$ in $C^2(\Gw\sms Q)$ and $\tl G_n(\cdot,x_0)\to\tl
G(\cdot,x_0)$ in $\mathrm{Lip}\,(\Gw)$. If $\gz_n=Z\tl G_n$ we have,

$$\BAL-\int_{\Gw_n}v\Gd\gz_n\,dx&=\int_{\bdw_n} v\prt_{\bmn}\gz\, dS=\int_{\bdw_n} vZ\prt_{\bmn}\tl G_n(\gx,x_0)\,dS\\
&=\int_{\bdw_n} vZP^{\Gw_n}(x_0,\gx)\, dS= \int_{\bdw_n} vZ\,d\gw_n.
\EAL$$
By \eqref{tr-w-conv},
$$\int_{\bdw_n} vZ\,d\gw_n\to\int_{\bdw}Z\,d\mu.$$
On the other hand, in view of \eqref{ibp0}, we have
$$\Gd \gz_n=\tl G_n\Gd Z +Z\Gd\tl G_n + 2\nabla Z\cdot \nabla \tl G_n\to \Gd Z$$
in $L^1_\gr(\Gw)$; therefore,
$$-\int_{\Gw_n}v\Gd\gz_n\,dx\to -\int_{\Gw}v\Gd\gz\,dx.$$
\qed

We denote by $\GTM_q=\GTM_q(\bdw)$ the set of q-good measures  on the boundary. A positive solution $u$ of \eqref{M1} in $\Gw$  possesses a boundary trace $\mu\in \GTM(\bdw)$ if and only if
\begin{equation}\label{INT}
 \myint{\Gw}{}u^q\gr dx<\infty
  \end{equation}
 \cite[Proposition 4.1]{MV7}. In this case $\mu\in \GTM_q$.
 \smallskip

The following statements can be proved in the same way as in the
case of smooth domains. For the proof in that case see \cite{MV2}.

\medskip
\nind{\bf I.} $\GTM_q(\bdw)$ is a linear space and
$$\mu\in \GTM_q(\bdw)\iff  |\mu|\in \GTM_q(\bdw).$$

\nind{\bf II.} If $\set{\mu_n}$ is an increasing \seq of measures in $\GTM_q(\bdw)$ and $\mu:=\lim \mu_n$ is a finite measure then $\mu\in \GTM_q(\bdw)$.

\bprop{edge condition} Let $\mu$ be a  bounded measure on $\prt P$.
($\mu$ may be a signed measure.)  For $i=1,\dots,N,\;j=1,\dots,
n_i$, we define the measure $\mu_{k,j}$ on $d_{A_{k,j}}$ by,
$$\mu_{k,j}=\mu\txt{on ${L_{k,j}}$,}\q \mu_{k,j}=0 \txt{on $d_{A_{k,j}}\sms {L_{k,j}}.$}$$
Then $\mu\in \GTM_q(\prt P)$, i.e., problem
\begin{equation}\label{Pbvp}
  -\Gd u+u^q=0 \txt{in $P$,} u=\mu \txt{on $\prt P$}
\end{equation}
possesses a solution, if and only if, $\mu_{k,j}$ is a q-good
measure relative to $D_{A_{k,j}}$ for all $(k,j)$ as above. \es

\Proof In view of statement {\bf I} above, it is sufficient to prove
the proposition in the case that $\mu$ is non-negative. This is
assumed hereafter. If $\mu\in \GTM_q(\prt P)$ then any measure $\nu$ on
$\prt P$ \sth $0\leq \nu\leq \mu$ is a q-good measure relative to
$P$. Therefore
$$\mu\in \GTM_q(\prt P) \Lra \mu'_{k,j}:=\mu\chi\indx{L_{k,j}}\in \GTM_q(\prt P).$$
 Assume that $\mu\in \GTM_q(\prt P)$ and let $u_{k,j}$ be the solution of \eqref{Pbvp} when $\mu$ is replaced by $\mu'_{k,j}$. Denote by $u'_{k,j}$  the extension of $u_{k,j}$ by zero to the k-dihedron $D_{A_{k,j}}$. Then $u'_{k,j}$ is a subsolution of \eqref{M1} in
$D_{A_{k,j}}$ with boundary data $\mu_{k,j}$. In the present case
there always exists a supersolution, e.g. the maximal solution of
\eqref{M1} in $D_{A_{k,j}}$ vanishing outside $d_{A_{k,j}}\sms\bar
L_{k,j}$. Therefore there exists a solution $v_{k,j}$ of this
equation in $D_{A_{k,j}}$ with boundary data $\mu_{k,j}$, i.e.,
$\mu_{k,j}$ is q-good relative to $D_{A_{k,j}}$.

Next assume that $\mu\in \GTM(\prt P)$ and that $\mu_{k,j}$ is
q-good relative to $D_{A_{k,j}}$ for every $(k,j)$ as above. Let
$v_{k,j}$ be the solution of \eqref{M1} in $D_{A_{k,j}}$ with
boundary data $\mu_{k,j}$. Then $v_{k,j}$ is a supersolution of
problem \eqref{Pbvp} with $\mu$ replaced by $\mu'_{k,j}$ and \consy
there exists a solution $u_{k,j}$ of this problem. It follows that
$$w:=\max\set{u_{k,j}: k=1,\dots,N,\;j=1,\dots, n_k}$$
is a subsolution  while
$$\bar w:=\sum_{k=1,\dots,N,\;j=1,\dots, n_k}u_{k,j}$$
is a supersolution of \eqref{Pbvp}. \Consy there exists a solution of this problem, i.e., $\mu\in \GTM_q(\prt P).$
\qed

\subsection{Removable singular sets and 'good measures', I}

We first introduce some standard elements associated to the Bessel capacities which are the natural way to characterize good measures or removable sets. For $\ga\in\BBR$, we denote by $G_\ga$ the Bessel kernel of order $\ga$, defined by
\begin{equation}\label{Bess1}
G_\ga(\xi)=\CF^{-1}\left((1+|\,.\,|^2)^{-\frac{\ga}{2}}\right)(\xi),
\end{equation}
where $\CF$ is the Fourier transform in the space $\CS'(\BBR^\ell)$ of moderate distributions in $\BBR^\ell$. For $1\leq p\leq\infty$, the Bessel space $L_{\ga,p}(\BBR^\ell)$ is defined by
\begin{equation}\label{Bess2}
L_{\ga,p}(\BBR^\ell)=\{f:f=G_\ga\ast g,: g\in L^p(\BBR^\ell)\},
\end{equation}
with norm
$$\norm f_{L_{\ga,p}}=\norm g_{L_{p}}=\norm {G_{-\ga}\ast f}_{L_{p}}.
$$

For $\ga,\gb\in\BBR$ and $1<p<\infty$, the mapping $f\mapsto G_\gb\ast f$ is an isomorphism from $L_{\ga,p}(\BBR^\ell)$ into
$L_{\ga+\gb,p}(\BBR^\ell)$.
Finally the Bessel spaces are connected to Besov and Sobolev spaces: when $\ga>0$ and $1<p<\infty$, it is known that if $\ga\in\BBN$, $L_{\ga,p}(\BBR^\ell)=W^{\ga,p}(\BBR^\ell)$ and if $\ga\notin\BBN$, then $L_{\ga,p}(\BBR^\ell)=B^{\ga,p}(\BBR^\ell)$, with equivalent norms (see e.g. \cite{Ca}, \cite{St}).

The Bessel capacity $C_{\ga,p}^{\BBR^\ell}$ ($\ga>0$, $p\geq 1$) is defined by the following rules: if $K\subset\BBR^\ell$ is compact
\begin{equation}\label{Bess3}
C_{\ga,p}^{\BBR^\ell}(K)=\inf\{\norm f^p_{L_{\ga,p}}:f\in\CS(\BBR^\ell), f\geq \chi_{K}\}.
\end{equation}
If $G$ is open
\begin{equation}\label{Bess4}
C_{\ga,p}^{\BBR^\ell}(G)=\sup\{C_{\ga,p}^{\BBR^\ell}(K):K\subset G,\,K \text{ compact}\}.
\end{equation}
If $A$ is any set
\begin{equation}\label{Bess5}
C_{\ga,p}^{\BBR^\ell}(A)=\inf\{C_{\ga,p}^{\BBR^\ell}(G):A\subset G,\,G \text{ open}\}.
\end{equation}

Note that the capacity of any non-empty set is positive if and only if $\ga>\frac{\ell}{p}$ because of Sobolev-Besov imbedding theorem.

\bprop{admi-conv} Let $A$ be a \Lip domain on $S^{k-1}$, $2\leq
k\leq N-1$, and let $D_A$ be the k-dihedron with opening $A$.
 Let $\gm\in\GTM(\prt D_{A})$ be a positive measure with compact support contained in $d_{A}$
 (= the edge of $D_A$). Assume that $\mu$ is q-good relative to $D_A$. Let $R>1$ be large enough so that
$\supp\mu\sbs B^{N-k}_R(0)$  and let $u$ be the solution of
\eqref{M1} in $D_A^R$ with trace $\mu$ on $d_A^R$ and trace zero
on $\prt D_A^R\sms d_A^R$. Then:

 \medskip
\nind{\rm (i)} For every non-negative $\eta\in
C_0^\infty(B_{3R/4}^{N-k}(0))$, 
\begin{equation}\label{inverse-ineq}\BAL
&\left( \int_{d_{A}^R}\eta^{q'} d\gm \right)\leq cM^{q'}
\int_{D_{A}^R}{}u^q\gr dx+\\
&\q+ cM^{q'}\left(\int_{D_{A}^R}{}u^q\gr dx\right)^{\rec{q}}
\left(1+M^{-1}\norm{\eta}_{L^{q'}(d_{A}^R)}\right). \EAL
\end{equation}
where $M=\norm{\eta}\indx{L^\infty}$ and $\gr$ is the first
eigenfunction of $-\Gd$ in $D_A^R$ normalized by $\gr(x_0)=1$ at
some point $x_0\in D_A^R$. The constant $c$ depends only on $N,q,k,
x_0,\gl_1,R$ where $\gl_1$ is the first eigenvalue.

\smallskip
 \nind{\rm (ii)} For any compact set
 $E\subset d_{A}$,

\begin{equation}\label{BVP-A2}
C^{N-k}_{s,q}(E)=0\Longrightarrow \mu (E)=0, \q s=2-\frac{\gk_{+}+k}{q'},
\end{equation}
where $C^{N-k}_{s,q}$ denotes the Bessel capacity with the indicated
indices in $\BBR^{N-k}$. \es

\Remark If we replace $D_A^R$ by $D_A\cap B_{\tl
R}^k(0)\cap B_{ R}^{N-k}(0)$,  $\tl R>1$,  then the constant $c$ in
(i) depends on $\tl R$ but {\em not on $ R$.}

\medskip
\Proof  We identify $d_A$ with $\BBR^{N-k}$ and  use the notation
$$x=(x',x'')\in \BBR^k\ti\BBR^{N-k},\q  y=|x'|.$$

 Let $\eta\in C_{0}^\infty(\BBR^{N-k})$ and let $R$ be large enough so that
$\supp\eta\sbs B^{N-k}_{R/2}(0)$.
 Let $w=w_R(t,x'')$   be the solution of the following problem in $\BBR_+\ti B^{N-k}_{R}(0)$:
  \begin{equation}\label{heat1}\BAL
 \prt_{t}w-\Gd_{x''} w&=0&&\text{in }\BBR^+\ti B_R^{N-k}(0),\\
 w(0,x'')&=\eta (x'')&&\text{in }B_R^{N-k},\\
 w(t,x'')&=0&&\text{on }\prt B^{N-k}_{R}(0).
 \EAL\end{equation}
Thus $w_R(t,\cdot)=S_R(t)[\eta]$ where $S_R(t)$ is the semi-group operator corresponding to the above problem. Denote,
   \begin{equation}\label{heat2}
    H_R[\eta](x',x'')=w_R(|x'|^2,x'')=S_R(y^2)[\eta](x''),\q y:=|x'|.
   \end{equation}

We assume, as we may, that $R>1$. Let $\gr^R$ be the first eigenfunction of $-\Gd_{x''}$ in the ball $B^{N-k}_{R}(0)$ normalized by $\gr^R(0)=1$ and let $\gr_A$ be the first eigenfunction of $-\Gd_{x'}$ in $C_{A}$
(where $C_{A}$ denotes the cone with opening $A$ in $\BBR^k$) normalized so that $\gr_A(x'_0)=1$ at some point $x'_0\in S_A$. Then $\gr^R\gr_{A}$ is the first eigenfunction of $-\Gd$ in $\set{x\in D_{A}:|x''|<R}$. Note that $\gr^R\leq 1$ and $\gr^R\to 1$  as $R\tin$ in $C^2(I)$ for any bounded set $I\sbs \BBR^{N-k}$.

Let $h\in C^\infty(\BBR)$ be a monotone decreasing function \sth $h(t)=1$ for $t<1/2$ and $h(t)=0$ for $t>3/4$. Put
$$\psi_R(x')=h(|x'|/R)$$
and
\begin{equation}\label{gz1}
   \gz_R:=\gr_{A}\psi_R H_R[\eta]^{q'}.
\end{equation}
If $\gr_A^R$ is the first eigenfunction (normalized at $x_0$) of $D_A^R:=D_A\cap \Gg_R$ ($\Gg_R$
as in \eqref{GgR}) then
\begin{equation}\label{rhoAR}
  \gr_A\psi_R\leq c\gr_A^R
\end{equation}
 and $\gr^R\gr_A^R$ is the first eigenfunction in $D_A^R$.

Hereafter we shall drop the index $R$  in $\gz_R, H_R, w_R$ but keep it in the other notations in order
 to avoid confusion.

We shall verify that $\gz\in D_A^R$. To this purpose we compute,
  \begin{equation}\label{heat2'}\BAL
\Gd\gz=&-\gl_{1}(\gr_A\psi_R) H[\eta]^{q'}+(\gr_A\psi_R)\Gd H[\eta]^{q'}+2\nabla(\gr_A\psi_R)\cdot\nabla H[\eta]^{q'}\\
=&-\gl_{1}\gz+q'(\gr_A\psi_R) (H[\eta])^{q'-1}\Gd H[\eta]\\
&+q(q'-1)(\gr_A\psi_R) (H[\eta])^{q'-2}| \nabla H[\eta]|^2\\
&+2q'(H[\eta])^{q'-1}\nabla(\gr_A\psi_R)\cdot\nabla H[\eta].
 \EAL\end{equation}
In addition,
$$\BA {l}\nabla H[\eta]=\nabla_{x'} H[\eta]+\nabla_{x''} H[\eta]
 =\prt_{y}H[\eta]\myfrac{x'}{y}+\nabla_{x''} H[\eta]\\[2mm]\phantom{\nabla H[\eta]}
 =2y\prt _{t}w(y^2,x'')\myfrac{x'}{y}+\nabla_{x''} H[\eta](x',x'')
 \EA$$
and \consy (recall that $y$ stands for $|x'|$),
$$\BAL &\nabla H[\eta]\cdot\nabla(\gr_A\psi_R)\\
&=2\prt_{t}w(y^2,x'')x'\cdot\Big (\psi_R\big(|x'|^{\kappa_+-1}(\kappa_+\frac{x'}{y}\gw_k(x'/y)+ |x'|\nabla \gw_k(x'/y))\big)+\gr_A\nabla \psi_R\Big)\\
&=2\kappa_+\prt_{t}w(y^2,x'')
|x'|^{\kappa_+}\gw_k(x'/y)=2\prt_{t}w(y^2,x'')(\kappa_+\gr_A \psi_R+
\gr_A x'\cdot\nabla\psi_R). \EAL$$ Since $w=w_R$ vanishes for
$|x''|=R$ and $\eta=0$ in a neighborhood of this sphere,
$|\prt_{t}w(y^2,x'')|\leq c\gr^R$. As $\psi_R$ vanishes for
$|x'|>3R/4$ we have $\gr_A  \nabla\psi_R\leq c\gr_A^R$. Therefore
$$
\abs{\nabla H[\eta]\cdot\nabla\gr_A}\leq c \gr^R\gr_A^R
$$
and, in view of \eqref{heat2'},
  \begin{equation}\label{Gdz}
  \abs{\Gd\gz}\leq c \gr^R\gr_A^R.
\end{equation}
  Thus $\gz\in X(D_{A}^R)$ and \consy
  \begin{equation}\label{int-eq-8.1}
\int_{D_{A}^R}\left(-u\Gd\gz+u^q\gz\right)dx=-\int_{D_{A}^R}\BBK[\gm]\Gd\gz dx.
\end{equation}
Since $q(q'-1)\gr_A (H[\eta])^{q'-2}| \nabla H[\eta]|^2\geq 0$, we have
  \begin{equation}\label{heat3}\BAL
 & \abs{\int_{D_{A}^R}u\Gd\gz dx}\\
&\leq  \int_{D_{A}^R}u\left(\gl_{1}\gz+q' (H[\eta])^{q'-1}\left(\gr|\Gd H[\eta]|+2|\nabla\gr.\nabla H[\eta]| \right)\right) dx\\
&\leq  \int_{D_{A}^R}{}u\left(\gl_{1}\gz + q' \gz^{1/q}\left(\gr^{1/q'}|\Gd H[\eta]|+2\gr^{-1/q}|\nabla\gr.\nabla H[\eta]| \right)\right) dx\\
&\leq  \left(\int_{D_{A}^R}{}u^q\gz dx\right)^{\rec{q}}
\left(\gl_{1}\left(\myint{D_{A}^R}{}\gz dx\right)^{\rec{q'}}+q'\norm{L[\eta]}_{L^{q'}(D_{A}^R)}\right) \EAL\end{equation}
where
  \begin{equation}\label{heat4}
  L[\eta]=\gr^{1/q'}|\Gd H[\eta]|+2\gr^{-1/q}|\nabla\gr.\nabla H[\eta]|.
  \end{equation}
By \rprop{IBP}
    \begin{equation}\label{heat4'}-\myint{D_{A}^R}{}\BBK[\gm]\Gd \gz dx=\myint{d_{A}^R}{}\eta^{q'}d\gm.
 \end{equation}
  Therefore
    \begin{equation}\label{heat5}\BAL
&\left( \int_{d_{A}^R}\eta^{q'} d\gm \right)\leq
\int_{D_{A}^R}{}u^q\gz dx+\\
&\q+ \left(\int_{D_{A}^R}{}u^q\gz dx\right)^{\rec{q}}
\left(\gl_{1}\left(\myint{D_{A}^R}{}\gz
dx\right)^{\rec{q'}}+q'\norm{L[\eta]}_{L^{q'}(D_{A}^R)}\right). \EAL
\end{equation}

 Next we prove that
\begin{equation}\label{heat6}
\norm{L[\eta]}_{L^{q'}(D_{A}^R)}\leq
C\norm\eta_{W^{s,q'}(\BBR^{N-k})}
\end{equation}
starting with the estimate of the first term on the right hand side
of \eqref{heat4}.

 $$\BA {l}\Gd H[\eta]=\Gd_{x'} H[\eta]+\Gd_{x''} H[\eta]
 =\prt^2_{y}H[\eta]+\myfrac{k-1}{y}\prt_{y}H[\eta]+\Gd_{x''} H[\eta]\\[2mm]\phantom{\nabla H[\eta]}
 =2y^2\prt_{tt}w(y^2,x'')+k\prt_{t}w(y^2,x'')+\Gd_{x''} H[\eta]\\[2mm]\phantom{\nabla H[\eta]}
 =2y^2\prt_{tt}w(y^2,x'')+(k+1)\prt_{t}w(y^2,x'').
 \EA$$
 Then
  $$\BA {l}\myint{\BBR^N}{}\gr \abs{\Gd H[\eta]}^{q'}dx\leq
 c\myint{0}{1}\myint{\BBR^{N-k}}{} \abs{\prt_{tt}w(y^2,x'')}^{q'}dx''y^{\gk_{+}+2q'+k-1}dy\\[4mm]
 \phantom{\myint{\BBR^N}{}\gr \abs{\Gd H[\eta]}^{q'}dx}
 +
 c\myint{0}{1}\myint{\BBR^{N-k}}{} \abs{\prt_{t}w(y^2,x'')}^{q'}dx''y^{\gk_{+}+k-1}dy\\[4mm]
 \phantom{\myint{\BBR^N}{}\gr \abs{\Gd H[\eta]}^{q'}dx}
 \leq c\myint{0}{1}\myint{\BBR^{N-k}}{} \abs{\prt_{tt}w(t,x'')}^{q'}dx''t^{(\gk_{+}+k)/2+q'}\myfrac{dt}{t}
 \\[4mm]
 \phantom{\myint{\BBR^N}{}\gr \abs{\Gd H[\eta]}^{q'}dx}
 +
c\myint{0}{1}\myint{\BBR^{N-k}}{} \abs{\prt_{t}w(t,x'')}^{q'}dx''t^{(\gk_{+}+k)/2}\myfrac{dt}{t}
\\[4mm]
 \phantom{\myint{\BBR^N}{}\gr \abs{\Gd H[\eta]}^{q'}dx}
\leq c\myint{0}{1}\norm{t^{2-(1-\frac{\gk_{+}+k}{2q'}))}\myfrac{d^2S(t)[\eta]}{dt^2}}^{q'}_{L^{q'}(\BBR^{N-k})}\myfrac{dt}{t}\\[4mm]
 \phantom{\myint{\BBR^N}{}\gr \abs{\Gd H[\eta]}^{q'}dx}
+ c\myint{0}{1}\norm{t^{1-(1-\frac{\gk_{+}+k}{2q'})}\myfrac{dS(t)[\eta]}{dt}}^{q'}_{L^{q'}(\BBR^{N-k})}\myfrac{dt}{t}.
 \EA$$
Put $\gb=\frac{\gk_{+}+k}{2q'}$ and note that
$0<\gb=\rec{2}(2-s)<1$. By standard interpolation theory,
$$\BAL
&\myint{0}{1}\norm{t^{1-(1-\gb)}\myfrac{dS(t)[\eta]}{dt}}^{q'}_{L^{q'}(\BBR^{N-k})}\myfrac{dt}{t}\\
&\approx\norm {\eta}^{q'}_{\left[W^{2,q'},L^{q'}\right]_{1-\gb,q'}}\approx
\norm\eta^{q'}_{W^{2(1-\gb),q'}
(\BBR^{N-k})},
\EAL$$
and
$$\BAL
&\myint{0}{1}\norm{t^{2-(1-\gb))}\myfrac{d^2S(t)[\eta]}{dt^2}}^{q'}_{L^{q'}(\BBR^{N-k})}\myfrac{dt}{t}\\
&\approx\norm
{\eta}^{q'}_{\left[W^{4,q'},L^{q'}\right]_{\frac{1}{2}(1-\gb),q'}}\approx
\norm\eta^{q'}_{W^{2(1-\gb),q'} (\BBR^{N-k})}. \EAL$$

The second term on the right hand side of \eqref{heat4} is estimated
in a similar way:
 $$\BAL &\int_{\BBR^N}\gr^{-q'/q}\abs{\nabla H[\eta]\cdot\nabla\gr}^{q'}dx\leq c\myint{0}{1}\int_{\BBR^{N-k}}
 \abs{\prt _{t}w(y^2,x'')}^{q'}dx'y^{\gk_{+}+k-1}dy\\
  &\leq c\myint{0}{1}\int_{\BBR^{N-k}}
 \abs{\prt _{t}w(t,x'')}^{q'}dx't^{\frac{\gk_{+}+k}{2}}\myfrac{dt}{t}\\
&\leq c\myint{0}{1}\norm{t^{1-(\frac{1}{2}-\gb)}\myfrac{d S(t)[\eta]}{dt}}_{L^{q'}(\BBR^{N-k})}^{q'}\myfrac{dt}{t}\\
&\approx\norm{\eta}^{q'}_{W^{2(1-\gb),q'}(\BBR^{N-k})}.
\EAL$$

This proves  \eqref{heat6}. Further, \eqref{heat5} and \eqref{heat6}
imply \eqref{inverse-ineq}.

We turn to the proof of part (ii). Let  $E$ be a closed subset of
$B^{N-k}_{R/2}(0)$ such that $C^{N-k}_{s,q'}(E)=0$. Then there
exists a \seq $\{\eta_{n}\}$ in $C^\infty_{0}(d_{A})$ such that
$0\leq\eta_{n}\leq 1$, $\eta_{n}=1$ in a neighborhood of $E$ (which
may depend on $n$), $\supp\eta_n\sbs B^{N-k}_{3R/4}(0)$ and
$\norm{\eta_{n}}_{W^{s,q'}}\to 0$. Then, by \eqref{heat6},
$$\norm{L[\eta_{n}]}_{L^{q'}(D_{A}^R)}\to 0.$$
Furthermore
$$\norm{w}\indx{L^{q'}((0,R)\ti B_R^{N-k}(0))}\leq c\norm{\eta_n}\indx{L^{q'}(B_R^{N-k}(0))}$$
 and \consy
 $$H[\eta_n]\to 0
 \txt{in}  L^{q'}(D_A^R).$$
 (Here we use the fact that $k\geq2$.) In addition
 $$0\leq H[\eta_n]\leq
 1,\q H[\eta_n]\leq c(R-|x'|)$$
 with a constant $c$ independent of $n$. Hence (see
 \eqref{rhoAR})
 $$\gz_{n,R}:=\gr_{A}\psi_R H[\eta_n]^{q'}\leq \gr^R\gr_{A}\psi_R H[\eta_n]^{q'-1}\leq
 \gr^R\gr_{A}^R H[\eta_n]^{q'-1}.$$
As $u^q\gr^R\gr_A^R\in L^1(D_A^R)$ we obtain,
$$\lim_{n\to\infty}\myint{D_{A}}{}u^q\gz_{n}dx=0.
$$
This fact and \eqref{heat5} imply that
$$\int_{d_{A}^R}\eta_n^{q'} d\gm\to 0.$$
As $\eta_n=1$ on a \ngh of $E$ in $\BBR^{N-k}$ it follows that
$\mu(E)=0$. \qed

\bprop{edge-bvp} Let $D_A$ be a k-dihedron, $1\leq k< N$. Let $k_+$
be as in \eqref{kappa2} and let $q^*_c$ and $q_c$ be as in
\rprop{admp0} and \rprop{admpk} respectively. Assume that $q_c\leq
q<q^*_c$. A measure $\mu\in \GTM(\prt D_A)$, with compact support
contained in $d_A$, is q-good relative to $D_A$ if and only if $\mu$
vanishes on every Borel set $E\sbs d_A$ \sth $C_{s,q'}(E)=0$, where
$s=2-\frac{k+\kappa_+}{q'}$. \es

\Remark We shall use the notation $\mu\prec C_{s,q'}$ to say that
$\mu$ vanishes on any Borel set $E\sbs(d_A)$ \sth $C_{s,q'}(E)=0$.

In the case $k=N$: $D_A=C_A$ (= the cone with vertex $0$ and opening
$A$ in $\BBR^k$) and $q_c=q^*_c$. By \cite{MV7} (specifically the results quoted in subsection 2.2)
$q_c=1-\frac{2}{\kappa_-}=\frac{N+\kappa_+}{N+\kappa_+ -2}$ and if $1<q<q_c$
then there exist solutions for every measure
$\mu=k\gd_P$, $P\in d_A$.

In the case $k=1$, $q^*_c=\infty$, $\kappa_+=1$ and
$q_c=\frac{N+1}{N-1}$. Thus $s=2/q$ and the statement of the theorem
is well known (see \cite{MV3}).

\medskip
\Proof In view of the last remark, it remains to deal only with
$2\leq k\leq N-1$. We shall identify $d_A$ with $\BBR^{N-k}$.

It is sufficient to prove the result for positive measures because
$\mu\prec C_{s,q'}$ if and only if $|\mu|\prec C_{s,q'}$. In
addition, if $|\mu|$ is a q-good measure then $\mu$ is a q-good
measure.

First we show that if $\mu$ is non-negative and q-good then
$\mu\prec C_{s,q'}$. If $E$ is a Borel subset of $\bdw$ then
$\mu\chi\indx{E}$ is q-good. If $E$ is compact and $C_{s,q'}(E)=0$
then, by \rprop{admi-conv}, $E$ is a removable set. This means that
the only positive solution of \eqref{M1} in $D_A$ \sth $\mu(\bdw\sms E)=0$ is the
zero solution. This implies that $\mu\chi\indx{E}=0$, i.e.,
$\mu(E)=0$. If $C_{s,q'}(E)=0$ but $E$ is not compact then
$\mu(E')=0$ for every compact set $E'\sbs E$. Therefore, we conclude
again that $\mu(E)=0$.

Next, assume that $\mu$ is a positive measure in $\GTM(\prt D_A)$
supported in a compact subset of $\BBR^{N-k}$.

If $\mu\in B^{-s,q}(\BBR^{N-k})$ then, by \rth{main1}, $\mu$ is admissible relative to $D_A\cap \Gg_{k,R}$,
for every $R>0$. (As before $\Gg_{k,R}$ is the cylinder with radius
$R$ around the 'axis' $\BBR^{N-k}$.) This implies that $\mu$ is
q-good relative to $D_A$.

If $\mu\prec C_{s,q'}$ then, by a theorem of Feyel and de la
Pradelle \cite{FD} (see also \cite{BP}), there exists a \seq
$\{\mu_n\}\sbs (B^{-s,q}(\BBR^{N-k}))_+$ \sth $\mu_n\uparrow\mu$. As
$\mu_k$ is q-good,  it follows that  $\mu$ is q-good. \qed

\bth{poly-good} Let $P$ be an $N$-dimensional polyhedron as
described in \rprop{edge condition}.  Let $\mu$ be a bounded measure
on $\prt P$, \(may be a signed measure\). Let
$k=1,\dots,N,\;j=1,\dots, n_k$, and let ${L_{k,j}}$ and $A_{k,j}$ be
defined as at the beginning of this section. Further, put
\begin{equation}\label{skj}
 s(k,j)=2-\frac{k+(\kappa_+)_{k,j}}{q'},
\end{equation}
where $(\kappa_+)_{k,j}$ is defined as in \eqref{kappa2} with
$A=A_{k,j}$.
Then $\mu\in \GTM_q(\prt P)$, i.e., $\mu$ is a good measure for
\eqref{M1} relative to $P$, if and only if, for every pair $(k,j)$
as above and every Borel set $E\sbs L_{k,j}$:\\
 If $1\leq k<N$ then
\begin{equation}\label{poly-good1}\BAL
 (q_c)_{k,j}\leq q<(q^*_c)_{k,j},\; C^{N-k}_{s(k,j),q'}(E)=0 &\Lra
 \mu(E)=0 \\
 q\geq (q_c^*)_{k,j} &\Lra
 \mu(L_{N,j})=0
\EAL\end{equation}
and if $k=N$, i.e., $L$ is a vertex,
\begin{equation}\label{poly-good2}
q\geq
(q_c)_{k,j}=\frac{N+2+\sqrt{(N-2)^2+4\gl_A}}{N-2+\sqrt{(N-2)^2+4\gl_A}}\Lra
\mu(L)=0.
\end{equation}
Here $(q_c^*)_{k,j}$ and $(q_c)_{k,j}$ are defined  as in \eqref{qk}
and \eqref{q-critk}respectively,  with $A=A_{k,j}$.

If $1<q<(q_c)_{k,j}$ then there is no restriction on
$\mu\chi\indx{L_{k,j}}$. \es

\Proof This is an immediate consequence of \rprop{edge condition}
and \rprop{edge-bvp} (see also the Remark following it). In the case
$k=N$, $L_{N,j}$ is a vertex and the condition says merely that for
$q\geq q_c(L_{N,j}$, $\mu$ does not charge the vertex. \qed

\subsection{Removable singular sets II}

\bprop{uqro} Let $A$ be a \Lip domain on $S^{k-1}$, $2\leq k\leq
N-1$, and let $D_A$ be the k-dihedron with opening $A$. Let $u$ be a
positive solution of \eqref{M1} in $D_A^R$, for some $R>0$.
Suppose that $F=\CS(u)\sbs d_A^R$ and let $Q$ be an open \ngh of $F$
\sth $\bar Q\sbs d_A^R$. \(Recall that $d_A^R=d_A\cap B^{N-k}_R(0)$
is an open subset of $d_A$.\) Let $\mu$ be the trace of $u$ on
$\CR(u)$.

Let $\eta\in W_0^{s,q'}(d_A^R)$ \sth
\begin{equation}\label{eta-cond}
 0\leq\eta \leq 1,\q \eta=0 \txt{on $Q$.}
\end{equation}
Employing the notation in the proof of \rprop{admi-conv}, put

\begin{equation}\label{gz-eta}
  \gz:=\gr_{A}\psi_RH_R[\eta]^{q'}.
\end{equation}
 Then
\begin{equation}\label{uqro1}
  \int_{D_A^R} u^q \gz\,dx\leq
  c(1+\norm{\eta}\indx{W^{s,q'}(d_A)})^{q'}+\mu(d_A^R\sms Q)^q,
\end{equation}
$c$ independent of $u$ and $\eta$.

\es

\Proof First we prove \eqref{uqro1} for $\eta\in C_0^\infty(d_A^R)$.
Let $\gs_0$ be a point in $A$ and let $\set{A_n}$ be a \Lip
exhaustion of $A$. 
If $0<\ge<\dist(\prt A,\prt A_n)=\bar \ge_n$ then

$$\ge\gs_0+C_{A_n}\sbs C_A.$$
Denote
$$D_A^{R', R''}=D_A\cap[|x'|<R']\cap[|x''|<R''].$$
Pick a \seq $\{\ge_n\}$ decreasing to zero \sth
$0<\ge_n<\min(\bar\ge_n/2^n, R/8)$. Let $u_n$ be the function given
by
$$u_n(x'x'')=u(x'+\ge_n\gs_0,x'') \forevery x\in D_{A_n}^{R_n,R},
\q R_n=R-\ge_n.$$ Then $u_n$ is a solution of \eqref{M1} in
$D_{A_n}^{R_n,R}$ belonging to $C^2(\bar D_{A_n}^{R_n,R})$ and we
denote its boundary trace by $h_n$. Let
$$\gz_n:=\gr_{A_n}\psi_RH_R[\eta]^{q'},$$
with $\psi_R$ and $H_R[\eta]$ as in the proof of \rprop{admi-conv}.
By \rprop{IBP}
    \begin{equation}\label{heat4''}-\int_{D_{A_n}^{R_n,R}}\BBP[h_n]\Gd \gz_n dx=\int_{B_R^{N-k}(0)}\eta^{q'}h_nd\gw_n
 \end{equation}
 where $\gw_n$ is the harmonic measure on $d_{A_n}^R$ relative to $D_{A_n}^{R_n,R}$. (Note that $d_{A_n}^R=d_{A}^R$
 and we may identify it with  $B_R^{N-k}(0)$.) Hence
  \begin{equation}\label{eq-8.2}
\int_{D_{A_n}^{R_n,R}}\left(-u_n\Gd\gz_n+u_n^q\gz_n\right)dx=-\int_{B_R^{N-k}(0)}\eta^{q'}h_n
\,d\gw_n.
\end{equation}
Further,
$$\int_{B_R^{N-k}(0)}\eta^{q'}h_n \,d\gw_n\to
\int_{B_R^{N-k}(0)}\eta^{q'}d\mu\leq \mu(d_A^R\sms Q),$$ because
$\eta=0$ in $Q$. By \eqref{heat3}, \eqref{heat6} we obtain,
  \begin{equation}\label{ineq-8.2}\BAL
&\abs{\int_{D_{A_n}^{R_n,R}}u_n\Gd\gz_n\,dx}\leq\\
c\Big(&\int_{D_{A_n}^{R_n,R}}u_n^q\gz_n dx\Big)^{\rec{q}}
\Big(\Big(\int_{D_{A_n}^{R_n,R}}\gz_n
dx\Big)^{\rec{q'}}+\norm{\eta}_{W^{s,q'}(B^{N-k}_R(0))}\Big).
\EAL\end{equation}

\nind From the definition of $\gz_n$ it follows that
$$\int_{D_{A_n}^{R_n,R}}\gz_n\,dx\leq \int_{D_{A_n}^{R_n,R}}\gr_n\,dx\to \int_{D_{A}^{R}}\gr\,dx,$$
where $\gr$ (resp. $\gr_n$)  is the first eigenfunction of $-\Gd$ in
$D_A^R$ \(resp. $D_{A_n}^{R_n,R}$\) normalized by $1$ at some
$x_0\in D_{A_1}^{R_1,R}$. Therefore, by \eqref{eq-8.2},

$$\BAL&\int_{D_{A_n}^{R_n,R}}u_n^q\gz_ndx\leq c\Big(&\int_{D_{A_n}^{R_n,R}}u_n^q\gz_n dx\Big)^{\rec{q}}
\big(1+\norm{\eta}_{W^{s,q'}(B^{N-k}_R(0))}\big) + \mu(d_A^R\sms
Q).\EAL$$
This implies
\begin{equation}\label{ineq-8.3'}
   \int_{D_{A_n}^{R_n,R}}u_n^q\gz_ndx\leq
   c\big(1+\norm{\eta}_{W^{s,q'}(B^{N-k}_R(0))}\big)^{q'}+\mu(d_A^R\sms
   Q)^q.
\end{equation}
To verify this fact, put
$$\BAL m=\Big(\int_{D_{A_n}^{R_n,R}}u_n^q\gz_ndx\Big)^{1/q},\; b=\mu(d_A^R\sms
Q),\; a=c\big(1+\norm{\eta}_{W^{s,q'}(B^{N-k}_R(0))}\big)\EAL$$ so
that \eqref{ineq-8.3'} becomes
$$m^q-am-b\leq 0.$$
If $b\leq m$ then
$$m^{q-1}-a-1\leq 0.$$
Therefore,
$$m\leq (a+1)^{\rec{q-1}}+b$$
which implies \eqref{ineq-8.3'}. Finally, by  the lemma of Fatou we
obtain \eqref{uqro1} for $\eta\in C_0^\infty$. By continuity we
obtain the inequality for any $\eta\in W^{s,q'}_0$ satisfying
\eqref{eta-cond}. \qed

\bth{removable} Let $A$ be a \Lip domain on $S^{k-1}$, $2\leq k\leq
N-1$, and let $D_A$ be the k-dihedron with opening $A$.
  Let $E$ be a compact subset of $d_A^R$ and let $u$ be a non-negative solution
 of \eqref{M1} in $D_A^R$ (for some $R>0$) \sth $u$ vanishes on
 $\prt D_A^R\sms E$. Then
\begin{equation}\label{REM1}
C^{N-k}_{s,q'}(E)=0,\q s=2-\frac{\gk_{+}+k}{q'} \Longrightarrow u=0,
\end{equation}
where $C^{N-k}_{s,q'}$ denotes the Bessel capacity with the
indicated indices in $\BBR^{N-k}$. \es

\Proof By \rprop{admi-conv}, \eqref{REM1} holds under the additional
assumption
\begin{equation}\label{REM2'}
\int_{D_A^R}u^q\gr_R\gr_A^R dx<\infty.
\end{equation}
 Indeed, by \cite[Proposition 4.1]{MV7},
\eqref{REM2'} implies that the solution $u$ possesses a boundary
trace $\mu$ on $\prt D_A^R$. By assumption, $\mu(\prt D_A^R\sms
E)=0$. Therefore, by \rprop{edge-bvp}, the fact that
$C^{N-k}_{s,q'}(E)=0$ implies that $\mu(E)=0$. Thus $\mu=0$ and
hence $u=0$.

We show that, under the conditions of the theorem, if
$C^{N-k}_{s,q'}(E)=0$ then \eqref{REM2'} holds.

By \rprop{uqro}, for every $\eta\in W_0^{s,q'}(d_A^R)$ \sth $0\leq
\eta\leq 1$ and $\eta=0 \txt{in a \ngh of $E$,}$
\begin{equation}\label{ineq-8.3}
   \int_{D_{A}^{R}}u^q\gz\,dx\leq
   c\big(1+\norm{\eta}_{W^{s,q'}(B^{N-k}_R(0))}\big)^{q'},
\end{equation}
for $\gz$ as in \eqref{gz-eta}. (Here we use the assumption that
$u=0$ on $\prt D_A^R\sms E$.)

 Let $a>0$ be sufficiently small so that $E\sbs
B_{(1-4a)R}^{N-k}(0)$. Pick a \seq $\{\phi_n\}$  in
$C_0^\infty(\BBR^{N-k})$ \sth, for each $n$, there exists a \ngh
$Q_n$ of $E$,
   $\bar Q_n\sbs B^{N-k}_{(1-3a)R}(0)$ and
\begin{equation}\label{phin-1}\BAL
   &0\leq \phi_n\leq 1 \txt{everywhere,} \phi_n=1 \txt{in $Q_n$,}\\
   &\tl\phi_n:=\phi_n\chi\indx{[|x''|<(1-2a)R]}\in
   C_0^\infty(\BBR^{N-k}),\\
&\big\|\tl\phi_n\big\|_{W^{s,q'}(\BBR^{N-k})}\to 0\txt{as} n\tin\\
&\eta_n:=(1-\phi_n)\lfloor\indx{[|x''|<R]}\in
   C_0^\infty(d_A^R),\\ & \eta_n=0 \;\text{in  }[(1-a)R<|x''|<R].
 \EAL\end{equation}

\nind  Such a \seq exists because $C^{N-k}_{s,q'}(E)=0$. Applying
\eqref{ineq-8.3} to $\eta_n$ we obtain,
\begin{equation}\label{ineq-8.4}
 \sup  \int_{D_{A}^{R}}u^q\gz_n\,dx\leq c<\infty,
\end{equation}
where $\gz_n=\gr_A\psi_RH_R^{q'}[\eta_n]$ (see \eqref{gz-eta}). By
taking a \sseq we may assume that $\{\eta_n\}$ converges (say to
$\eta$) in $L^{q'}(B_{R}^{N-k}(0))$ and \consy $H[\eta_n]\to
H[\eta]$ in the sense that
$$H_R[\eta_n](x',\cdot)=w_{n,R}(y^2,\cdot)\to
w_{R}(y^2,\cdot)=H_R[\eta](x',\cdot) \txt{in $L^{q'}$}$$ uniformly
\wrto $y=|x'|$.  It follows that
\begin{equation}\label{ineq-8.5}
\int_{D_{A}^{R}}u^q\gz\,dx<\infty,\q  \gz=\gr_A\psi_RH_R^{q'}[\eta].
\end{equation}

As $\tl\phi_n\to 0$ in $W^{s,q'}(\BBR^{N-k})$ it follows that
$\phi_n\to 0$ and hence $\eta_n\to 1$ a.e. in
$B_{(1-2a)R}^{N-k}(0)$. Thus $\eta=1$ in this ball, $\eta=0$ in
$[(1-a)R<|x''|<R]$ and $0\leq \eta\leq 1$ everywhere.

 \Consy, given
$\gd>0$, there exists an $N$-dimensional \ngh $O$ of $d_A\cap
B_{(1-2a)R}^{N-k}(0)$ \sth
$$1-\gd<H_R[\eta]<1 \txt{and} 1-\gd<\psi_R/\gr_A^R<1 \txt{in}O.$$
Therefore \eqref{ineq-8.5} implies that
\begin{equation}\label{ineq-8.6}
   \int_{D_{A}^{(1-3a)R}}u^q\gr^R\gr_A^R\,dx\leq c<\infty.
\end{equation}
Recall that the trace of $u$ on  $\prt D_A^R\sms d_A^{(1-4a)R}$ is
zero. Therefore $u$ is bounded in $D_A^R\sms D_{A}^{(1-3a)R}$. This
fact and \eqref{ineq-8.6} imply \eqref{REM2'}. \qed

\bdef{ro-capacity} Let $\Gw$ be a bounded \Lip domain. Denote by
$\gr$ the first eigenfunction of $-\Gd$ in $\Gw$ normalized by
$\gr(x_0)=1$ for a fixed point $x_0\in\Gw$.

For every compact set $K\sbs \bdw$ we define
$$ M_{\gr, q}(K)=\set{\mu\in \GTM(\bdw):\mu\geq 0,\;\mu(\bdw\sms
K)=0,\; \BBK[\mu]\in L^q_\gr(\Gw)}$$
and

$$\tl C_{\gr,q'}(K)=\sup\set{\mu(K)^q:\;\mu\in M_{\gr, q}(K),\;\int_\Gw
\BBK[\mu]^{q}\gr\,dx=1}.$$

Finally we denote by $C_{\gr,q'}$ the outer measure generated by the
above functional. \es

The following statement is verified by standard arguments:

\blemma{cgrq} For every compact $K\sbs \bdw$, $C_{\gr,q'}(K)=\tl
C_{\gr,q'}(K)$. Thus $C_{\gr,q'}$ is a capacity and,
\begin{equation}\label{cgrq}
    C_{\gr,q'}(K)=0 \iff M_{\gr, q}(K)=\{0\}.
\end{equation}
\es

\bth{removable set} Let $\Gw$ be a bounded polyhedron in $\BBR^N$. A
compact set $K\sbs \bdw$ is removable if and only if
\begin{equation}\label{remov1}
  C_{s(k,j),q'}(K\cap L_{k,j})=0,
\end{equation}
for $k=1,\cdot,N$ $j=1,\cdots, n_k$, where $s(k,j)$ is defined as in
\eqref{skj}. This condition is equivalent to
\begin{equation}\label{remov1'}
  C_{\gr,q'}(K)=0.
\end{equation}
A measure $\mu\in\GTM(\bdw)$ is q-good if and only if it does not
charge sets with $C_{\gr,q'}$-capacity zero.
 \es

\Proof The first assertion is an immediate consequence of
\rprop{edge condition} and \rth{removable}. The second assertion
follows from the fact that
$$C_{\gr,q'}(K\cap L_{k,j})=C_{s(k,j),q'}(K\cap L_{k,j}).$$
The third assertion follows from \rth{poly-good} and the previous
statement. \qed




\end {document}